\documentclass{article}
\usepackage{amsmath}
\usepackage{amssymb}
\begin{document}
\def\e#1\e{\begin{equation}#1\end{equation}}
\def\ea#1\ea{\begin{align}#1\end{align}}
\def\eq#1{{\rm(\ref{#1})}}
\newtheorem{thm}{Theorem}[section]
\newtheorem{prop}[thm]{Proposition}
\newtheorem{cor}[thm]{Corollary}
\newenvironment{dfn}{\medskip\refstepcounter{thm}
\noindent{\bf Definition \thesection.\arabic{thm}\ }}{\medskip}
\newenvironment{proof}[1][,]{\medskip\ifcat,#1
\noindent{\it Proof.\ }\else\noindent{\it Proof of #1.\ }\fi}
{\hfill$\square$\medskip}
\def\dim{\mathop{\rm dim}}
\def\mod{\mathop{\rm mod}}
\def\Re{\mathop{\rm Re}}
\def\Im{\mathop{\rm Im}}
\def\vol{\mathop{\rm vol}}
\def\Tr{\mathop{\rm Tr}}
\def\GL{\mathop{\rm GL}}
\def\U{\mathbin{\rm U}}
\def\SL{\mathop{\rm SL}}
\def\SO{\mathop{\rm SO}}
\def\SU{\mathop{\rm SU}}
\def\sn{{\textstyle\mathop{\rm sn}}}
\def\cn{{\textstyle\mathop{\rm cn}}}
\def\dn{{\textstyle\mathop{\rm dn}}}
\def\sech{{\textstyle\mathop{\rm sech}}}
\def\cosech{\mathop{\rm cosech}}
\def\ge{\geqslant} 
\def\le{\leqslant} 
\def\cal{\mathcal}
\def\bb{\mathbb}   
\def\R{\mathbin{\mathbb R}}
\def\Z{\mathbin{\mathbb Z}}
\def\N{\mathbin{\mathbb N}}
\def\Q{\mathbin{\mathbb Q}}
\def\C{\mathbin{\mathbb C}}
\def\CP{\mathbb{CP}}
\def\RP{\mathbb{RP}}
\def\u{\mathfrak{u}} 
\def\su{\mathfrak{su}} 
\def\sl{\mathfrak{sl}}
\def\gl{\mathfrak{gl}}
\def\al{\alpha}
\def\be{\beta}
\def\ga{\gamma}
\def\de{\delta}
\def\ep{\epsilon}
\def\th{\theta}
\def\ze{\zeta}
\def\la{\lambda}
\def\ka{\kappa}
\def\vp{\varphi}
\def\si{\sigma}
\def\La{\Lambda}
\def\Om{\Omega}
\def\Si{\Sigma}
\def\om{\omega}
\def\d{{\rm d}}
\def\pd{\partial}
\def\db{{\bar\partial}}
\def\ts{\textstyle}
\def\w{\wedge}
\def\br{\buildrel}
\def\sm{\setminus}
\def\ov{\overline}
\def\ot{\otimes}
\def\iy{\infty}
\def\ra{\rightarrow}
\def\longra{\longrightarrow}
\def\t{\times}
\def\ha{{\textstyle\frac{1}{2}}}
\def\op{\oplus}
\def\ti{\tilde}
\def\bs{\boldsymbol}
\def\ms#1{\vert#1\vert^2}
\def\bms#1{\bigl\vert#1\bigr\vert^2}
\def\md#1{\vert #1 \vert}
\def\bmd#1{\bigl\vert #1 \bigr\vert}
\def\an#1{\langle#1\rangle}
\def\ban#1{\bigl\langle#1\bigr\rangle}
\title{Special Lagrangian 3-folds and integrable systems}
\author{Dominic Joyce \\ Lincoln College, Oxford}
\date{}
\maketitle

\section{Introduction}
\label{is1}

This is the sixth in a series of papers \cite{Joyc1,Joyc2,Joyc3,Joyc4,Joyc5}
constructing explicit examples of special Lagrangian submanifolds
(SL $m$-folds) in $\C^m$. The principal motivation for the series is 
to study the singularities of SL $m$-folds, especially when $m=3$. 
This paper also has a second objective, which is to connect SL $m$-folds 
with the theory of integrable systems, and to arouse interest in 
special Lagrangian geometry within the integrable systems community. 

We begin in \S\ref{is2} with a brief introduction to {\it special
Lagrangian submanifolds} in $\C^m$, which are a class of real 
$m$-dimensional minimal submanifolds in $\C^m$, defined using
calibrated geometry. Section \ref{is3} then gives a rather longer 
introduction to {\it harmonic maps} $\psi:S\ra\CP^{m-1}$, where $S$ 
is a Riemann surface. Such maps form an {\it integrable system}, and
have a complex and highly-developed theory involving the Toda lattice 
equations, loop groups, and classification using spectral curves.

Section \ref{is4} explains the connection of this with special 
Lagrangian geometry. Let $N$ be a {\it special Lagrangian cone} 
in $\C^3$, and set $\Si=N\cap{\cal S}^5$. Then $\Si$ is a 
{\it minimal Legendrian surface} in ${\cal S}^5$, and so the 
image of a {\it conformal harmonic map} $\phi:S\ra{\cal S}^5$ 
from a Riemann surface $S$. The projection $\psi=\pi\circ\phi$ 
of $\phi$ from ${\cal S}^5$ to $\CP^2$ is also conformal and 
harmonic, with Lagrangian image.

Thus, $\psi$ can be analyzed in the integrable systems framework 
of \S\ref{is3}. As the image of $\psi$ is Lagrangian there is a 
simplification, in which the $\SU(3)$ Toda lattice equation reduces 
to the {\it Tzitz\'eica equation}, and the spectral curve acquires 
an extra symmetry. We use the integrable systems theory to give 
{\it parameter counts} for the expected families of SL $T^2$-cones 
in~$\C^3$.

In \S\ref{is5} we give an explicit construction of special 
Lagrangian cones $N$ in $\C^3$, involving two commuting o.d.e.s, 
and reducing to constructions given in \cite{Joyc1,Joyc2} in 
special cases. Taking the intersection with ${\cal S}^5$, we 
obtain families of explicit conformal harmonic maps 
$\phi:\R^2\ra{\cal S}^5$ and $\psi:\R^2\ra\CP^2$. Under some 
circumstances we can solve the conditions for these maps to 
be {\it doubly-periodic} in $\R^2$, and so to push down to 
harmonic maps $T^2\ra{\cal S}^5$ and~$T^2\ra\CP^2$.

Section \ref{is6} analyzes this family of harmonic maps 
$\psi:\R^2\ra\CP^2$ from the point of view of integrable systems.
We find that for generic initial data $\psi$ is {\it superconformal},
and explicitly determine its harmonic sequence, Toda and Tzitz\'eica
solutions, algebra of polynomial Killing fields, and spectral curve. 
In \S\ref{is7} we generalize the ideas of \S\ref{is5} to give a new 
construction of special Lagrangian 3-folds in $\C^3$, which involves 
{\it three} commuting o.d.e.s, and reduces to the construction of 
\S\ref{is5} in a special case.

We end with an open problem. The SL 3-folds of \S\ref{is7} look very 
similar to those of \S\ref{is5}, and share many of the hallmarks of 
integrable systems -- commuting o.d.e.s, elliptic functions, conserved 
quantities. The author wonders whether these examples can also be 
explained in terms of some higher-dimensional integrable system, and 
indeed whether the special Lagrangian equations themselves are in 
some sense integrable.
\medskip

\noindent{\it Acknowledgements.} I would like to thank Ian McIntosh,
Karen Uhlenbeck, Robert Bryant and Mark Haskins for helpful 
conversations. I would also like to thank the organizers of the 
`Integrable Systems in Differential Geometry' conference, Tokyo,
2000, where I began to understand the material of \S\ref{is3} 
and~\S\ref{is6}.

\section{Special Lagrangian submanifolds in $\C^m$}
\label{is2}

We begin by defining {\it calibrations} and {\it calibrated 
submanifolds}, following Harvey and Lawson~\cite{HaLa}.

\begin{dfn} Let $(M,g)$ be a Riemannian manifold. An {\it oriented
tangent $k$-plane} $V$ on $M$ is a vector subspace $V$ of
some tangent space $T_xM$ to $M$ with $\dim V=k$, equipped
with an orientation. If $V$ is an oriented tangent $k$-plane
on $M$ then $g\vert_V$ is a Euclidean metric on $V$, so 
combining $g\vert_V$ with the orientation on $V$ gives a 
natural {\it volume form} $\vol_V$ on $V$, which is a 
$k$-form on~$V$.

Now let $\vp$ be a closed $k$-form on $M$. We say that
$\vp$ is a {\it calibration} on $M$ if for every oriented
$k$-plane $V$ on $M$ we have $\vp\vert_V\le \vol_V$. Here
$\vp\vert_V=\al\cdot\vol_V$ for some $\al\in\R$, and 
$\vp\vert_V\le\vol_V$ if $\al\le 1$. Let $N$ be an 
oriented submanifold of $M$ with dimension $k$. Then 
each tangent space $T_xN$ for $x\in N$ is an oriented
tangent $k$-plane. We say that $N$ is a {\it calibrated 
submanifold} if $\vp\vert_{T_xN}=\vol_{T_xN}$ for all~$x\in N$.
\label{is2def1}
\end{dfn}

It is easy to show that calibrated submanifolds are automatically
{\it minimal submanifolds} \cite[Th.~II.4.2]{HaLa}. Here is the 
definition of special Lagrangian submanifolds in $\C^m$, taken
from~\cite[\S III]{HaLa}.

\begin{dfn} Let $\C^m$ have complex coordinates $(z_1,\dots,z_m)$, 
and define a metric $g$, a real 2-form $\om$ and a complex $m$-form 
$\Om$ on $\C^m$ by
\e
\begin{split}
g=\ms{\d z_1}+\cdots+\ms{\d z_m},\quad
\om&=\frac{i}{2}(\d z_1\w\d\bar z_1+\cdots+\d z_m\w\d\bar z_m),\\
\text{and}\quad\Om&=\d z_1\w\cdots\w\d z_m.
\end{split}
\label{is2eq}
\e
Then $\Re\Om$ and $\Im\Om$ are real $m$-forms on $\C^m$. Let
$L$ be an oriented real submanifold of $\C^m$ of real dimension 
$m$, and let $\th\in[0,2\pi)$. We say that $L$ is a {\it special 
Lagrangian submanifold} of $\C^m$ if $L$ is calibrated with 
respect to $\Re\Om$, in the sense of Definition \ref{is2def1}.
We will often abbreviate `special Lagrangian' by `SL', and 
`$m$-dimensional submanifold' by `$m$-fold', so that we shall
talk about SL $m$-folds in~$\C^m$. 
\end{dfn}

As in \cite{Joyc1} there is also a more general definition 
of special Lagrangian submanifolds involving a {\it phase} 
${\rm e}^{i\th}$, but we will not use it in this paper. Harvey 
and Lawson \cite[Cor.~III.1.11]{HaLa} give the following 
alternative characterization of special Lagrangian submanifolds.

\begin{prop} Let\/ $L$ be a real $m$-dimensional submanifold 
of $\C^m$. Then $L$ admits an orientation making it into an
SL submanifold of\/ $\C^m$ if and only if\/ $\om\vert_L\equiv 0$ 
and\/~$\Im\Om\vert_L\equiv 0$.
\label{is2prop}
\end{prop}

Note that an $m$-dimensional submanifold $L$ in $\C^m$ is 
called {\it Lagrangian} if $\om\vert_L\equiv 0$. Thus special 
Lagrangian submanifolds are Lagrangian submanifolds satisfying 
the extra condition that $\Im\Om\vert_L\equiv 0$, which is how 
they get their name.

\section{Harmonic maps and integrable systems}
\label{is3}

A map $\phi:S\ra M$ of Riemannian manifolds is {\it harmonic} if 
it extremizes the energy functional $\int_S\ms{\d\phi}\d V$. When 
$S$ is 2-dimensional, the energy is conformally invariant, so that 
we may take $S$ to be a Riemann surface. In this case, if $\phi$ is 
conformal, then $\phi$ is harmonic if and only if $\phi(S)$ is 
{\it minimal}\/ in $M$. Thus, harmonic maps are closely connected
to minimal surfaces.

We shall describe a relationship, due to Bolton, Pedit and Woodward 
\cite{BPW}, between a special class of harmonic maps $\psi:S\ra\CP^{m-1}$
called {\it superconformal}\/ harmonic maps, and solutions of the 
{\it Toda lattice equations} for $\SU(m)$. Then we will explain
how superconformal maps can be studied using {\it loop groups} and
{\it loop algebras}. 

This leads to the definition of {\it polynomial Killing fields} and a 
special class of superconformal maps called {\it finite type}, which
include all maps from $T^2$. Finally we explain how to associate a
Riemann surface called the {\it spectral curve} to each finite type
superconformal map, and that finite type harmonic maps can be classified
in terms of algebro-geometric data including the spectral curve.

This is a deep and complex subject, and we cannot do it justice in
a few pages. A good general reference on the following material is 
Fordy and Wood \cite{FoWo}, in particular, the articles by 
Bolton and Woodward \cite[p.~59--82]{FoWo}, McIntosh 
\cite[p.~205--220]{FoWo} and Burstall and 
Pedit~\cite[p.~221--272]{FoWo}.

\subsection{The harmonic sequence and superconformal maps}
\label{is31}

Suppose $S$ is a connected Riemann surface and $\psi:S\ra\CP^{m-1}$ 
a harmonic map. Then the {\it harmonic sequence} $(\psi_k)$ of 
$\psi$ is a sequence of harmonic maps $\psi_k:S\ra\CP^{m-1}$ with 
$\psi_0=\psi$, defined in Bolton and Woodward \cite[\S 1]{BoWo}.
Each $\psi_k:S\ra\CP^{m-1}$ defines a holomorphic line subbundle $L_k$ 
of the trivial vector bundle $S\t\C^m$, where a section $s$ of 
$L_k$ is defined to be holomorphic if $\pd s/\pd\bar z$ is orthogonal 
to~$L_k$. 

The $\psi_k$ and $L_k$ are characterized by the following property. 
If $U$ is an open subset of $S$ and $z$ a holomorphic coordinate on 
$U$, then any nonzero holomorphic section $\phi_0$ of $L_0$ over $U$
may be extended uniquely to a sequence of nonzero holomorphic sections
$\phi_k$ of $L_k$ over $U$ satisfying
\e
\begin{gathered}
\ban{\phi_k,\phi_{k+1}}=0,\quad
\frac{\pd\phi_k}{\pd z}=\phi_{k+1}+\frac{\pd}{\pd z}
\bigl(\log\ms{\phi_k}\bigr)\phi_k\\
\text{and}\quad
\frac{\pd\phi_k}{\pd\bar z}=-\frac{\ms{\phi_k}}{\ms{\phi_{k-1}}}\phi_{k-1}
\quad\text{for all $k$,}
\end{gathered}
\label{is3eq1}
\e
where $\an{\,,\,}$ is the standard Hermitian product on $\C^m$. 
(Actually one should allow the $\phi_k$ to be meromorphic, but we will 
ignore this point.) If $(\phi_k)$, $(\phi_k')$ both satisfy \eq{is3eq1} 
then $\phi_k'=f\phi_k$ for some holomorphic $f:U\ra\C^*$ and all $k$, 
where $\C^*=\C\sm\{0\}$. Thus $\psi_k=[\phi_k]$ is independent of the 
choice of~$\phi_0$.

Note that $\psi_k$ may not be defined for all $k\in\Z$. For if 
$\psi_k:S\ra\CP^{m-1}$ is holomorphic then $\frac{\pd\phi_k}{\pd z}=0$, 
so that $\phi_{k+1}=0$ and $\psi_{k+1}$ is undefined, and the sequence 
terminates above at $\psi_k$. Similarly, if $\psi_k$ is antiholomorphic 
then $\psi_{k-1}$ is undefined, and the sequence terminates below 
at~$\psi_k$.

If $\psi_k$ exists for all $k\in\Z$ then $\psi$ is called 
{\it non-isotropic}. Otherwise $\psi$ is called {\it isotropic}. 
Isotropic maps $\psi:S\ra\CP^{m-1}$ were studied by Eells and Wood 
\cite{EeWo}. They all arise by projection from certain holomorphic 
maps into a complex flag manifold, and so are fairly easy to 
understand and construct.

Harmonic sequences have strong orthogonality properties. Two points 
in $\CP^{m-1}$ are called {\it orthogonal}\/ if the corresponding lines in 
$\C^m$ are orthogonal at all points, and two maps $\psi_j,\psi_k:
S\ra\CP^{m-1}$ are called orthogonal if $\psi_j(s)$ and $\psi_k(s)$ are 
orthogonal in $\CP^{m-1}$ for all $s\in S$. 

Bolton and Woodward \cite[Prop.~2.4]{BoWo} show that if some set of 
$l$ consecutive terms in a harmonic sequence $(\psi_k)$ are mutually 
orthogonal, then {\it every} set of $l$ consecutive terms of $(\psi_k)$ 
are mutually orthogonal. A harmonic map $\psi:S\ra\CP^{m-1}$ and its 
harmonic sequence $(\psi_k)$ are both called $l$-{\it orthogonal}\/ if 
every set of $l$ consecutive terms are mutually orthogonal. 

Clearly, every harmonic sequence is 2-orthogonal. It is easy to show 
that $\psi=\psi_0$ is conformal if and only if $\psi_1$ and $\psi_{-1}$ 
are orthogonal. Therefore, $\psi$ is 3-orthogonal if and only if it is 
conformal, and then all the elements $\psi_k$ of the harmonic sequence 
are also conformal.

The maximum number of mutually orthogonal elements of $\CP^{m-1}$
is $m$, and so a harmonic map $\psi:S\ra\CP^{m-1}$ is at most
$m$-orthogonal. A nonisotropic, $m$-orthogonal harmonic map 
$\psi:S\ra\CP^{m-1}$ is called {\it superconformal}. The harmonic 
sequence of a superconformal map $\psi$ is {\it periodic}, 
with period $m$, so that $\psi_{k+m}=\psi_k$ for all~$k$.

A nonisotropic, conformal harmonic map $\psi:S\ra\CP^2$ 
is superconformal, as $\psi$ is 3-orthogonal because it is 
conformal, from above. Thus, every conformal harmonic map
$\psi:S\ra\CP^2$ is either isotropic or superconformal.

\subsection{The Toda lattice equations}
\label{is32}

The {\it Toda lattice equations} for $\SU(m)$ may be
written as follows. For all $k\in\Z$, let $\chi_k:\C\ra(0,\iy)$
be differentiable functions satisfying
\begin{gather}
\chi_{k+m}=\chi_k \quad\text{for all $k\in\Z$,}\quad 
\chi_0\chi_1\cdots\chi_{m-1}\equiv 1,\quad\text{and}
\label{is3eq2}\\
\frac{\pd^2}{\pd z\pd\bar z}\bigl(\log\chi_k\bigr)
=\chi_{k+1}\chi_k^{-1}-\chi_k\chi_{k-1}^{-1}\quad\text{for all $k\in\Z$.}
\label{is3eq3}
\end{gather}
They are important integrable equations in mathematical physics, and 
large classes of solutions to them may be constructed using loop 
algebra methods. 

We shall show how to construct a solution of \eq{is3eq2}--\eq{is3eq3}
from a superconformal map $\psi:S\ra\CP^{m-1}$. Use the notation of 
\S\ref{is31}, and suppose $\psi$ is superconformal. Define functions 
$\chi_k:U\ra(0,\iy)$ by $\chi_k=\ms{\phi_k}$. Using the fact that 
$\pd^2\phi_k/\pd z\pd \bar z=\pd^2\phi_k/\pd\bar z\pd z$, one can 
show using \eq{is3eq1} that
\begin{equation*}
\frac{\pd^2}{\pd z\pd\bar z}\log\ms{\phi_k}=
\frac{\ms{\phi_{k+1}}}{\ms{\phi_k}}-\frac{\ms{\phi_k}}{\ms{\phi_{k-1}}}.
\end{equation*}
Thus the $\chi_k$ satisfy~\eq{is3eq3}. 

To make the $\chi_k$ satisfy \eq{is3eq2} as well, we need to choose 
the coordinate $z$ and lifts $\phi_k$ more carefully. As $\psi$ is 
superconformal, the $\psi_k$ are periodic with period $m$. We
shall arrange for the lifts $\phi_k$ also to be periodic with period 
$m$. Then $\chi_{k+m}=\chi_k$ for all $k$, the first equation of 
\eq{is3eq2}. This can be done by a suitable choice of holomorphic 
coordinate~$z$. 

It is not difficult to show that the $\phi_k$ satisfy $\phi_{k+m}=
\xi\phi_k$ for all $k$, where $\xi$ is a nonzero holomorphic function 
on $U$. If we change to a new holomorphic coordinate $z'$ on $U$, then 
$\xi$ is replaced by 
\e
\xi'=\Bigl(\frac{\pd z'}{\pd z}\Bigr)^{-m}\xi.
\label{is3eq4}
\e
Thus, by changing coordinates we can arrange that $\xi\equiv 1$, so 
that $\phi_{k+m}=\phi_k$ for all $k$, as we want. A holomorphic 
coordinate $z$ on an open subset $U$ of $S$ with this property is called
{\it special}\/ \cite[p.~65]{FoWo}. Such coordinates are unique up to 
addition of a constant, and multiplication by an $m^{\rm th}$
root of unity.

Suppose from now on that $z$ is special, so that $\phi_{k+m}=\phi_k$ 
for all $k$. It remains to show that we can choose the $\phi_k$ such
that the second equation of \eq{is3eq2} holds. Regard the $\phi_k$ as 
complex column vectors, so that $(\phi_0\,\phi_1\cdots\phi_{m-1})$ is 
a complex $m\t m$ matrix. Then the determinant $\det(\phi_0\,\phi_1
\cdots\phi_{m-1})$ is a nonzero holomorphic function on~$U$. 

From above, the $\phi_k$ are defined uniquely up to multiplication
by some holomorphic function $f:U\ra\C^*$. By multiplying the 
$\phi_k$ by a suitable $f$ we can arrange that
\e
\det(\phi_0\,\phi_1\cdots\phi_{m-1})\equiv 1. 
\label{is3eq5}
\e
This fixes the $\phi_k$ uniquely up to multiplication by an $m^{\rm th}$ 
root of unity. As $\psi$ is superconformal, $\phi_0,\ldots,\phi_{m-1}$ 
are complex orthogonal in $\C^m$. It follows that
\begin{equation*}
\chi_0\chi_1\cdots\chi_{m-1}\equiv\ms{\phi_0}\ms{\phi_1}\cdots\ms{\phi_{m-1}}
\equiv\bms{\det(\phi_0\,\phi_1\cdots\phi_{m-1})}\equiv 1,
\end{equation*}
so that the second equation of \eq{is3eq2} holds. Thus equations
\eq{is3eq2}--\eq{is3eq3} hold, and the $\chi_k$ satisfy the Toda 
lattice equations for~$\SU(m)$.

When $S$ is a torus $T^2$ and $\psi:S\ra\CP^{m-1}$ a superconformal
harmonic map, from \cite[Cor.~2.7]{BPW} and \cite[p.~67-8]{FoWo} 
there exists a global special holomorphic coordinate $z$ on the 
universal cover $\C$ of $T^2$, which then yields a solution 
$(\chi_k)$ of the Toda lattice equations \eq{is3eq2}--\eq{is3eq3} 
on the whole of $\C$. In particular, there are no `higher order 
singularities', and the $\phi_k$ and $\xi$ do not have zeros or poles.

\subsection{Toda frames and the reconstruction of $\psi$} 
\label{is33}

Above we saw how to construct a solution $(\chi_k)$ of the
Toda lattice equations for $\SU(m)$ out of a superconformal 
harmonic map $\psi:S\ra\CP^{m-1}$. We shall now explain how 
to go the other way, and reconstruct $\psi$ from $(\chi_k)$. 
We continue to use the same notation.

Define $F:U\ra\GL(m,\C)$ by $F=(f_0\,f_1\cdots f_{m-1})$, where 
$f_j=\md{\phi_j}^{-1}\phi_j$. Then as $\phi_0,\ldots,\phi_{m-1}$
are complex orthogonal and $\det(\phi_0\,\phi_1\cdots\phi_{m-1})=1$,
$F$ actually maps $U\ra\SU(m)$. We call $F$ a {\it Toda frame} 
for $\psi$ on $U$, \cite[p.~126]{BPW}. Define $\al$ to be the 
matrix-valued 1-form $F^{-1}\d F$ on $U$. Then by \eq{is3eq1} we 
find that $\al$ is given by
\e
\begin{pmatrix}
-\frac{i}{2}J\d\log\chi_0 & -\chi_1^{1/2}\!\chi_0^{-1/2}\!\d\bar z &&& 
\chi_0^{1/2}\!\chi_{m-1}^{-1/2}\d z \\
\chi_1^{1/2}\!\chi_0^{-1/2}\!\d z & -\frac{i}{2}J\d\log\chi_1 & \ddots \\
& \chi_2^{1/2}\!\chi_1^{-1/2}\!\d z & \ddots & 
-\chi_{m-2}^{1/2}\chi_{m-3}^{-1/2}\d\bar z \\
&& \ddots & -\frac{i}{2}J\d\log\chi_{m-2} & 
-\chi_{m-1}^{1/2}\chi_{m-2}^{-1/2}\d\bar z \\
-\chi_0^{1/2}\!\chi_{m-1}^{-1/2}\d\bar z &&& 
\chi_{m-1}^{1/2}\chi_{m-2}^{-1/2}\d z & 
-\frac{i}{2}J\d\log\chi_{m-1} 
\end{pmatrix},
\label{is3eq6}
\e
where $J$ is the complex structure on $U$. 

Now $\al$ is a 1-form on $U$ with values in $\su(m)$, so we 
may regard it as a {\it connection $1$-form} upon the trivial 
$\SU(m)$-bundle over $U$. The connection $\d+\al$ is automatically 
{\it flat}, as $\al$ is of the form $F^{-1}\d F$, so that $\al$
satisfies
\e
\d\al+\ha[\al\w\al]=0.
\label{is3eq7}
\e
Furthermore, $\al$ depends only on the solution $(\chi_k)$ of 
the Toda equations. 

Thus, to reconstruct $\psi$ from $(\chi_k)$, we proceed as follows.
Given $(\chi_k)$, we can write down the flat $\su(m)$-connection 
$\d+\al$ on $U$. Then we retrieve the Toda frame $F:U\ra\SU(m)$ by 
solving the equation $\d F=F\al$, which is in effect two commuting 
first-order linear o.d.e.s. If $U$ is simply-connected there exists 
a solution $F$, which is unique up to multiplication $F\mapsto AF$ 
by some~$A\in\SU(m)$. 

Finally we define $\psi=[f_0]$, where $f_0$ is the first column of $F$. 
In this way, any solution $(\chi_k)$ of the Toda lattice equations on 
a simply-connected open set $U$ in $\C$ generates a superconformal
map $\psi:U\ra\CP^{m-1}$, which is unique up to multiplication by
$A\in\SU(m)$, that is, up to automorphisms of~$\CP^{m-1}$.

\subsection{Loop groups and loops of flat connections}
\label{is34}

A large part of the integrable systems literature on harmonic maps
is formulated in terms of infinite-dimensional Lie groups known
as {\it loop groups}. If $G$ is a finite-dimensional Lie group, the 
{\it loop group} $LG$ is the group of smooth maps ${\cal S}^1\ra G$,
under pointwise multiplication and inverses, and the corresponding
{\it loop algebra} $L{\mathfrak g}$ is the Lie algebra of smooth maps
${\cal S}^1\ra{\mathfrak g}$, where $\mathfrak g$ is the Lie algebra of~$G$.

In the situation of \S\ref{is33}, for each $\la\in\C$ with $\md{\la}=1$, 
define an $\su(m)$-valued 1-form $\al_\la$ on $U$ to be
\e
\begin{pmatrix}
-\frac{i}{2}J\d\log\chi_0 & 
\!\!-\la\chi_1^{1/2}\!\chi_0^{-1/2}\!\d\bar z \!\!
&&& \!\!\la^{-1}\chi_0^{1/2}\!\chi_{m-1}^{-1/2}\d z \\
\la^{-1}\chi_1^{1/2}\!\chi_0^{-1/2}\!\d z\!\! 
& \!\!-\frac{i}{2}J\d\log\chi_1\!\! & \ddots \\
& \!\!\la^{-1}\chi_2^{1/2}\!\chi_1^{-1/2}\!\d z\!\! & \ddots & 
\!\!-\la\chi_{m-2}^{1/2}\chi_{m-3}^{-1/2}\d\bar z \\
&& \ddots & \!\!-\frac{i}{2}J\d\log\chi_{m-2}\!\! & 
\!\!-\la\chi_{m-1}^{1/2}\chi_{m-2}^{-1/2}\d\bar z \\
-\la\chi_0^{1/2}\!\chi_{m-1}^{-1/2}\d\bar z\!\! &&& 
\!\!\la^{-1}\chi_{m-1}^{1/2}\chi_{m-2}^{-1/2}\d z\!\!& 
\!\!-\frac{i}{2}J\d\log\chi_{m-1}
\end{pmatrix}.
\label{is3eq8}
\e
When $\la=1$ this coincides with the 1-form $\al$ of \eq{is3eq6}. 
Using the Toda lattice equations one can show that $\d+\al_\la$ is
also a flat $\SU(m)$-connection, so that
\e
\d\al_\la+\ha[\al_\la\w\al_\la]=0.
\label{is3eq9}
\e

Thus the family $\{\al_\la\}$ gives a {\it loop of flat connections}.
We can interpret this in loop group terms as follows. We defined the 
$\al_\la$ as an ${\cal S}^1$ family of 1-forms on $U\subseteq\C$ with 
values in $\su(m)$, but we can instead regard it as a single 1-form 
on $U$ with values in the loop algebra $L\su(m)$. So the $\al_\la$ 
give an $L\SU(m)$-connection on $U$, which turns out to be flat.

If $U$ is simply-connected, there exists a smooth 1-parameter family of 
maps $F_\la:U\ra\SU(m)$ with $F_\la^{-1}\d F_\la=\al_\la$, which
are unique up to multiplication $F_\la\mapsto A_\la F_\la$ by elements
$A_\la\in\SU(m)$. The family $\{F_\la\}$ is called an {\it extended 
Toda frame} for $\psi$. In loop group terms, we may interpret the 
$F_\la$ as a map $U\ra L\SU(m)$. It turns out that each $F_\la$ is the 
Toda frame of a superconformal harmonic map $\psi_\la:U\ra\CP^{m-1}$, 
where the special holomorphic coordinate on $U$ is $\la^{-1}z$ rather 
than~$z$.

Now if $\Phi:U\ra L\SU(m)$ is any smooth map, then 
$\d+\Phi^{-1}\d\Phi$ is a flat $L\SU(m)$-connection on $U$, or 
equivalently a loop of flat $\SU(m)$-connections on $U$. This 
gives an enormous family of loops of flat $\SU(m)$-connections 
on $U$, most of which have nothing to do with harmonic maps into 
$\CP^{m-1}$. The important thing about the family $\{\al_\la\}$ is 
that it has two special algebraic properties.

The first property is that we may write $\al_\la$ in the form
\e
\al_\la=(\al_1'\la+\al_0')\d z+(\al_{-1}''\la^{-1}+\al_0'')\d\bar z,
\label{is3eq10}
\e
where $\al_1',\al_0',\al_{-1}''$ and $\al_0''$ map $U\ra\su(m)^{\mathbb C}
=\sl(m,\C)$. This equation says two things. Firstly, as a Laurent
series in $\la$ we have $\al_\la=\al_1\la+\al_0+\al_{-1}\la^{-1}$. 
Secondly, if we decompose $\al_\la$ into $(1,0)$ and $(0,1)$ parts
as $\al_\la=\al_\la'\d z+\al_\la''\d\bar z$, then $\al_\la'=
\al_1'\la+\al_0'$, so that $\al_\la'$ has no $\la^{-1}$ component,
and $\al_\la''=\al_{-1}''\la^{-1}+\al_0''$, so that $\al_\la''$
has no $\la$ component.

The second property is this. Define $\ze={\rm e}^{2\pi i/m}$, so 
that $\ze^m=1$, and let $\Upsilon$ be the diagonal $m\t m$ matrix 
with entries $1,\ze^{-1},\ze^{-2},\ldots,\ze^{1-m}$. Then 
\e
\al_{\ze\la}=\Upsilon\al_\la\Upsilon^{-1}
\quad\text{for all $\la\in\C$ with $\md{\la}=1$.}
\label{is3eq11}
\e
That is, $\al_\la$ is equivariant under $\Z_m$-actions on ${\cal S}^1$ 
and~$\su(m)$.

It follows from Bolton, Pedit and Woodward \cite[\S 2]{BPW} that an 
${\cal S}^1$-family $\d+\al_\la$ of flat $\SU(m)$-connections on 
a simply-connected open subset $U\subseteq\C$ come from a solution
of the Toda lattice equations, and hence from a superconformal map 
$\psi:U\ra\CP^{m-1}$, if and only if the $\al_\la$ satisfy 
\eq{is3eq10}--\eq{is3eq11} and the additional condition that 
$\det(\al_1')$ is nonzero except at isolated points in~$U$.

\subsection{Polynomial Killing fields}
\label{is35}

Polynomial Killing fields were introduced by Ferus et al.\ 
\cite[\S 2]{FPPS} and used extensively by Burstall et al.\ 
\cite{BFPP}, but in a somewhat different situation to us. Our 
treatment is based on McIntosh \cite[App.~A]{McIn3} and Bolton, 
Pedit and Woodward~\cite[\S 3]{BPW}.

We shall work with the Lie algebra $\u(m)$ and its complexification
$\gl(m,\C)$ rather than $\su(m)$ and $\sl(m,\C)$. Let $\ze$ and $\Upsilon$ 
be as in \S\ref{is34}. For each $d\in\N$, let $\La_d\gl(m,\C)$ be the 
vector space of maps $\eta:\C^*\ra\gl(m,\C)$ of the form 
$\eta(\la)=\sum_{n=-d}^d\eta_n\la^n$, where $\eta_n\in\gl(m,\C)$,
which satisfy
\e
\eta(\ze\la)=\Upsilon\eta(\la)\Upsilon^{-1}
\quad\text{for all $\la\in\C^*$.}
\label{is3eq12}
\e
That is, $\eta$ is equivariant under the same $\Z_m$-actions as 
$\al_\la$ in~\eq{is3eq11}.

Let $\La_d\u(m)$ be the real vector subspace of $\eta\in\La_d\gl(m,\C)$ 
such that $\eta(\la)$ lies in $\u(m)$ for all $\la\in\C$ with $\md{\la}=1$. 
Then $\La_d\gl(m,\C)=\La_d\u(m)\ot_{\mathbb R}\C$. Note that by restricting 
$\eta$ to ${\cal S}^1$ in $\C^*$, we can regard $\La_d\gl(m,\C)$ and 
$\La_d\u(m)$ as finite-dimensional vector subspaces of the loop algebras 
$L\gl(m,\C)$ and~$L\u(m)$.

We define a {\it polynomial Killing field} on $U$ to be a map 
$\eta:U\ra\La_d\gl(m,\C)$ for some $d\in\N$ satisfying
\e
\d\eta=[\eta,\al_\la].
\label{is3eq13}
\e
We call $\eta$ {\it real}\/ if it maps to $\La_d\u(m)$ in 
$\La_d\gl(m,\C)$. We may write $\eta=\eta(\la,z)$ for $\eta\in\C^*$ 
and $z\in U$, and decompose $\eta$ as
\e
\eta(\la,z)=\sum_{n=-d}^d\eta_n(z)\la^n,
\quad\text{where $\eta_n$ maps $U\ra\gl(n,\C)$.}
\label{is3eq14}
\e

Using the decompositions \eq{is3eq10} and \eq{is3eq14} of $\al_\la$ and
$\eta$, it is easy to show that \eq{is3eq13} is equivalent to the equations
\ea
\frac{\pd\eta_n}{\pd z}&=[\eta_n,\al_0']+[\eta_{n-1},\al_1']
\quad\text{and}
\label{is3eq15}\\
\frac{\pd\eta_n}{\pd\bar z}&=[\eta_n,\al_0'']+[\eta_{n+1},\al_{-1}'']
\quad\text{for all $n$,}
\label{is3eq16}
\ea
where we set $\eta_n\equiv 0$ if~$\md{n}>d$. 

Define $\cal A$ to be the vector space of polynomial Killing fields. 
It is easy to see that the polynomial Killing fields form a Lie 
algebra under the obvious Lie bracket. In our case, where $\psi$ is 
superconformal, this Lie algebra is {\it abelian}, and the polynomial 
Killing fields form a {\it commutative algebra} under matrix 
multiplication \cite[p.~240-1]{McIn3}. The reason why we work with 
$\gl(m,\C)$ rather than $\sl(m,\C)$ is because $\sl(m,\C)$ is not 
closed under matrix multiplication.

Following \cite[p.~133]{BPW}, we say that $\psi$ is of {\it finite type} 
if there exists a real polynomial Killing field $\eta:U\ra\La_d\u(m)$ 
for some $d\equiv 1\mod m$ with $\eta_d=\al_1'$ and $\eta_{d-1}=2\al_0'$.
All finite type solutions may be obtained by integrating commuting
Hamiltonian o.d.e.s on the finite-dimensional manifold $\La_d\u(m)$,
and so are fairly well understood. By \cite[Cor.~3.7]{BPW}, every 
superconformal map corresponding to a doubly-periodic Toda solution 
on $\C$, and hence every superconformal $T^2$ in $\CP^{m-1}$, is of 
finite type. 

\subsection{Spectral curves}
\label{is36}

To each superconformal map $\psi:U\ra\CP^{m-1}$ of finite type one
can associate a Riemann surface known as a {\it spectral curve}.
There are two different definitions of spectral curve in use in
the literature. Here is the first. Fix $z\in U$, and define
\e
Y=\bigl\{(\la,[v])\in\C^*\t\CP^{m-1}:\text{$\forall\eta\in{\cal A}$, 
$\exists\mu\in\C$ with $\eta(\la,z)v=\mu v$}\bigr\}.
\label{is3eq17}
\e
The {\it spectral curve} as defined by Ferus, Pedit, Pinkall and 
Sterling \cite[\S 5]{FPPS} is the compactification $\ti Y$ of $Y$ 
in $\CP^1\t\CP^{m-1}$. In the generic case, $\ti Y$ is a compact,
nonsingular Riemann surface.

However, McIntosh \cite[App.~A]{McIn3} uses a different definition. As 
each $\eta\in{\cal A}$ satisfies \eq{is3eq12}, if $\bigl(\la,[v]\bigr)
\in Y$ then $\bigl(\ze\la,[\Upsilon v]\bigr)\in Y$. So define $\nu:Y\ra Y$ by
\e
\nu:\bigl(\la,[v]\bigr)\mapsto\bigl(\ze\la,[\Upsilon v]\bigr).
\label{is3eq18}
\e
Then $\nu^m=1$, and $\an{\nu}$ is a free $\Z_m$-action on $Y$, which
lifts to a free $\Z_m$-action on $\ti Y$. The spectral curve used by 
McIntosh is equivalent to $\ti X=\ti Y/\an{\nu}$. It is also 
generically a compact, nonsingular Riemann surface. 

The point of \eq{is3eq17} is that as $\cal A$ is a commutative 
algebra under matrix multiplication, the matrices $\eta(\la,z)$ 
can be simultaneously diagonalized for all $\eta\in{\cal A}$.
Assume that all the common eigenspaces are 1-dimensional. Then
$(\la,[v])\in Y$ if $[v]$ is an eigenspace of $\eta(\la,z)$ for 
all~$\eta\in{\cal A}$.

This suggests an alternative description of $Y$, using eigenvalues 
rather than eigenvectors. Pick $\eta\in{\cal A}$, and define
\e
Y'=\bigl\{(\la,\mu)\in\C^*\t\C:
\det\bigl(\mu I-\eta(\la,z)\bigr)=0\bigr\}.
\label{is3eq19}
\e
Define $\pi:Y\ra Y'$ by $\pi:\bigl(\la,[v]\bigr)\mapsto(\la,\mu)$ 
when $\eta(\la,z)v=\mu v$. Then $\pi$ is birational for generic 
$\eta$, and biholomorphic if $\eta$ generates $\cal A$ 
over~$\C[\la^mI,\la^{-m}I]$.

As a Riemann surface $\ti Y$ is independent of base-point $z\in U$,
but its embedding in $\CP^1\t\CP^{m-1}$ does depend on $z$. However,
$Y'$ is independent of $z$. This is because \eq{is3eq13} is 
equivalent to $\d(F_\la\eta F_\la^{-1})=0$, as $\al_\la=F_\la^{-1}
\d F_\la$. Thus $F_\la\eta F_\la^{-1}$ is independent of $z$, so
the eigenvalues of $\eta(\la,z)$ are independent of~$z$.

Spectral curves are used by McIntosh \cite{McIn1,McIn2,McIn3} to give 
a construction of all finite type harmonic maps $\psi:\C\ra\CP^{m-1}$,
and hence of all nonisotropic harmonic maps $\psi:T^2\ra\CP^{m-1}$.
Above we have only considered superconformal $\psi$, which are dealt 
with in \cite{McIn3}; the general case is studied in \cite{McIn1,McIn2},
and is rather more complicated.

By an explicit construction, McIntosh establishes a 1-1 correspondence 
between finite type nonisotropic harmonic maps $\psi:\C\ra\CP^{m-1}$ up 
to isometry, and quadruples of {\it spectral data} $(X,\si,\pi,{\cal L})$, 
where $X$ is a Riemann surface, $\si:X\ra X$ a real structure, 
$\pi:X\ra\CP^1$ a branched cover, and $\cal L$ a holomorphic line 
bundle over $X$, all satisfying certain conditions. Understanding 
the set of such $(X,\si,\pi,{\cal L})$ is a fairly straightforward 
problem in algebraic geometry. 

When $\psi$ is doubly-periodic it pushes down to $T^2$, and all 
nonisotropic harmonic maps $\psi:T^2\ra\CP^{m-1}$ arise in this way. 
This reduces the classification of nonisotropic harmonic tori in 
$\CP^{m-1}$ to the problem of understanding the double-periodicity 
conditions upon $(X,\si,\pi,{\cal L})$, which will be discussed
in~\S\ref{is43}.

\section{SL cones in $\C^3$ and the Tzitz\'eica equation}
\label{is4}

Suppose $N$ is a {\it special Lagrangian cone} in $\C^3$. Then 
$\Si=N\cap{\cal S}^5$ is a {\it minimal Legendrian surface} in 
${\cal S}^5$, and its projection $\pi(\Si)$ to $\CP^2$ is a 
{\it minimal Lagrangian surface} in $\CP^2$. Thus, $\Si$ and
$\pi(\Si)$ are the images of {\it conformal harmonic maps}
$\phi:S\ra{\cal S}^5$ and $\psi:S\ra\CP^2$, where $S$ is
a Riemann surface.

Therefore, we can apply the theory of \S\ref{is3} to $\psi$. It
turns out that as $\psi(S)$ is Lagrangian, the corresponding
solutions $\chi_k:U\ra(0,\iy)$ of the $\SU(3)$ Toda lattice 
equations simplify, coming from a single function $f:U\ra\R$ 
satisfying the {\it Tzitz\'eica equation}. This will be 
explained in~\S\ref{is41}.

The programme of \S\ref{is34}--\S\ref{is36} can then be applied,
to interpret finite type solutions of the Tzitz\'eica equation
in terms of {\it spectral data}. The difference with the $\SU(3)$
Toda lattice case is that the spectral curve $X$ acquires an 
extra symmetry, a holomorphic involution $\rho$ satisfying some 
compatibility conditions with the other data $\si,\pi,{\cal L}$.
The details can be found in Sharipov \cite{Shar} and Ma and Ma 
\cite{MaMa}, on which this section is based.

\subsection{Derivation of the Tzitz\'eica equation}
\label{is41}

Suppose $N$ is a special Lagrangian cone in $\C^3$. Define 
$\Si=N\cap{\cal S}^5$. Then $\Si$ is a {\it minimal Legendrian
surface} in ${\cal S}^5$, as $N$ is a minimal Lagrangian 3-fold 
in $\C^3$. It has a natural metric and orientation, and so 
inherits the structure of a Riemann surface. Let $S$ be $\Si$, 
regarded as an abstract Riemann surface, and $\phi:S\ra\Si$ 
the inclusion map.

Then $\phi$ is conformal by definition, and has minimal image, 
so it is harmonic. Let $\pi:{\cal S}^5\ra\CP^2$ be the Hopf 
projection, and define $\psi:S\ra\CP^2$ by $\psi=\pi\circ\phi$. 
As $\phi$ has Legendrian image it follows that $\psi$ is also 
a conformal harmonic map, with {\it Lagrangian image}. Therefore 
we can consider the harmonic sequence $(\psi_k)$ of $\psi$, as 
in~\S\ref{is31}. 

Let $z=x+iy$ be a holomorphic coordinate on an open subset 
$U\subset S$. Then $g\bigl(\phi,\frac{\pd\phi}{\pd x}\bigr)=
g\bigl(\phi,\frac{\pd\phi}{\pd y}\bigr)=0$ as $\md{\phi}\equiv 1$, 
and $\om\bigl(\phi,\frac{\pd\phi}{\pd x}\bigr)=\om\bigl(\phi,
\frac{\pd\phi}{\pd y}\bigr)=0$ as $N$ is Lagrangian. Hence, 
$\phi$ is complex orthogonal to $\frac{\pd\phi}{\pd\bar z}$, 
so $\phi$ is a {\it holomorphic} section of~$L_0$. 

Thus, by \S\ref{is31}, there is a unique sequence $(\phi_k)$ 
with $\phi_0=\phi$ satisfying equation \eq{is3eq1}. As 
$\md{\phi_0}\equiv 1$, we see that
\e
\phi_0=\phi,\quad
\phi_1=\frac{\pd\phi}{\pd z}\quad\text{and}\quad
\phi_{-1}=-\Big\vert\frac{\pd\phi}{\pd\bar z}\Big\vert^{-2}
\cdot\frac{\pd\phi}{\pd\bar z}.
\label{is4eq1}
\e
Suppose now that $\psi:S\ra\CP^2$ is {\it superconformal}, and that
$z$ is a {\it special} holomorphic coordinate. Then $\phi_k$ exists 
for all $k\in\Z$ and $\phi_{k+3}=\phi_k$, so that \eq{is4eq1}
determines $\phi_k$ for all $k$. Also, the functions $\chi_k=\ms{\phi_k}$
satisfy the {\it Toda lattice equations} \eq{is3eq2}--\eq{is3eq3} for~$m=3$.

Now in this case the Toda lattice equations simplify. For $\chi_0
\equiv 1$ as $\md{\phi_0}\equiv 1$, and so $\chi_2\equiv\chi_1^{-1}$ as 
$\chi_0\chi_1\chi_2\equiv 1$. Defining $f=\log\chi_1$ gives $\chi_{3k}=1$, 
$\chi_{3k+1}={\rm e}^f$ and $\chi_{3k-1}={\rm e}^{-f}$ for $f:U\ra\R$,
and then \eq{is3eq3} reduces to the single equation
\e
\frac{\pd^2f}{\pd z\pd\bar z}={\rm e}^{-2f}-{\rm e}^f.
\label{is4eq2}
\e

This is the elliptic version of the {\it Tzitz\'eica equation}. 
The corresponding hyperbolic equation first arose in 1910 in a 
study by Georges Tzitz\'eica \cite{Tzit} of a class of surfaces 
in $\RP^3$ now known as {\it affine spheres}. The equation was 
rediscovered in a solitonic context by Bullough and Dodd 
\cite{BuDo}, and amongst mathematical physicists is often 
known as the Bullough--Dodd equation.

We have shown that superconformal harmonic maps $\psi:S\ra\CP^2$ 
coming from special Lagrangian cones in $\C^3$ are related to
solutions of the Tzitz\'eica equation \eq{is4eq2}, in the same
way that general superconformal harmonic maps $\psi:S\ra\CP^{m-1}$ 
are related to solutions of the $\SU(m)$ Toda lattice equations.
The converse also applies, in that starting with a solution of
\eq{is4eq2} one can reconstruct a special Lagrangian cone in $\C^3$
using the Toda frame method of~\S\ref{is33}.

\subsection{Spectral data for the Tzitz\'eica equation}
\label{is42}

Suppose $f$ is a solution of the Tzitz\'eica equation. Then as 
in \S\ref{is41} we get a solution $(\chi_k)$ of the $\SU(3)$ 
Toda lattice equations. If $(\chi_k)$ is of finite type then 
as in \S\ref{is36} it has a set of {\it spectral data} 
$(X,\si,\pi,{\cal L})$. But because the $(\chi_k)$ come from a
solution of the Tzitz\'eica equation and so have a simplified
structure, the spectral data $(X,\si,\pi,{\cal L})$ has an
extra symmetry. 

It turns out that this symmetry is a {\it holomorphic involution} 
$\rho:X\ra X$. It commutes with $\si$, has the property that
if $\pi(x)=[1,\la]$ then $\pi\circ\rho(x)=[1,-\la]$ for $x\in X$
and $\la\in\C$, and lifts to a holomorphic involution of the line 
bundle~$\cal L$. 

Thus, just as there is a correspondence between finite type solutions
of the $\SU(m)$ Toda lattice equations and quadruples of spectral data
$(X,\si,\pi,{\cal L})$ satisfying certain conditions, there is also a 
correspondence between solutions of the Tzitz\'eica equation and
quintuples of spectral data $(X,\rho,\si,\pi,{\cal L})$ satisfying
certain conditions. 

This construction is used in two papers by Sharipov \cite{Shar} and 
Ma and Ma \cite{MaMa}. Sharipov considers `complex normal' surfaces in 
${\cal S}^5$, which in our terminology are just Legendrian surfaces. 
He shows that minimal Legendrian tori correspond to solutions of the 
Tzitz\'eica equation \eq{is4eq2}, gives the spectral data for finite 
type solutions, and sketches how to write the maps $\phi:\R^2\ra
{\cal S}^5$ in terms of Prym theta functions.

The paper by Ma and Ma is very similar. They consider `totally real'
surfaces in $\CP^2$, which in our terminology are just Lagrangian
surfaces. They show that minimal Lagrangian tori correspond to 
solutions of \eq{is4eq2}, give the spectral data, and write the 
maps $\psi:\R^2\ra\CP^2$ in terms of Prym theta functions. The
principal difference is that Ma and Ma give more proofs, more
detail, and more explicit formulae.

\subsection{Parameter counts for minimal tori in $\CP^2$}
\label{is43}

We shall now use the integrable systems set-up above to give 
parameter counts for the families of minimal tori in $\CP^2$ and 
of minimal Lagrangian tori in $\CP^2$ (equivalently, of special 
Lagrangian $T^2$-cones in $\C^3$ up to isomorphism). 

Similar parameter counts for case of minimal tori in $\CP^{m-1}$ are 
given by McIntosh in \cite[p.~516]{McIn2} and \cite[Th.~5]{McIn3},
and we follow his method, modifying it in the obvious way for the 
minimal Lagrangian case by requiring invariance under the holomorphic 
involution $\rho$. Each set of spectral data corresponds up to 
isometries of $\CP^2$ with a unique finite type harmonic map 
$\psi:\R^2\ra\CP^2$. We shall count the number of free parameters 
in the spectral data, and the number of restrictions for $\psi$ to 
be doubly-periodic.

First consider general minimal tori in $\CP^2$. The spectral data
for this is a quadruple $(X,\si,\pi,{\cal L})$. We suppose that
$X$ is a nonsingular Riemann surface of genus $p\ge 2$. (The
case $p\le 1$ is dealt with in \cite[\S 5]{McIn3}.) Now 
$\pi:X\ra\CP^1$ is a 3-fold branched cover, which we can regard 
as a {\it meromorphic function}, identifying $\CP^1$ with 
$\C\cup\{\iy\}$. It is required that $\pi^{-1}(0)$ and $\pi^{-1}(\iy)$ 
should be single points, and thus triple branch points. 

When $(X,\si,\pi,{\cal L})$ is generic, all other branch points
of $\pi$ will be double branch points, and there will be $2p$
of them by elementary topology. Let $\la_1,\ldots,\la_{2p}$
be the images of these branch points in $\C\sm\{0\}$. Then the 
$\la_k$ are distinct, and it is required \cite[Prop.~1]{McIn2}
that no $\la_k$ lies on the unit circle.

Now $\si$ acts on $X$, and if $\pi(x)=\la$ then $\pi\circ\si(x)=
1/\bar\la$. Clearly $\si$ must take double branch points to
double branch points, and so the set $\{\la_1,\ldots,\la_{2p}\}$
is closed under $\la\mapsto 1/\bar\la$. As no $\la_k$ lies on
the unit circle this swaps the $\la_k$ in pairs. Order the $\la_k$
so that $\md{\la_k}<1$ and $\la_{k+p}=1/\bar\la_k$ for~$k=1,\ldots,p$.

The triple $(X,\si,\pi)$ depends on $\la_1,\ldots,\la_p$ and discrete 
data. Thus there are $2p$ real parameters in $(X,\si,\pi)$. The
set of $\si$-invariant line bundles $\cal L$ has dimension $p$,
but $\cal L$ depends on a choice of base point in $\R^2$, so 
factoring out by translations in $\R^2$ shows that the choice
of $\cal L$ really represents $p-2$ degrees of freedom. 

From \cite[\S 5]{McIn3}, the double-periodicity conditions for 
$\psi:\R^2\ra\CP^2$ depend only on $(X,\si,\pi)$, being independent 
of $\cal L$. The condition for the Toda solution to be doubly-periodic 
is $2p-4$ rationality conditions, and for $\psi$ to be doubly periodic 
in $\CP^2$ is another 4 rationality conditions.

Effectively, this means that the moduli space ${\cal M}_p$ of data
$(X,\si,\pi)$ has dimension $2p$, and to be doubly-periodic requires 
that $2p$ real functions $f_1,\ldots,f_{2p}$ on this moduli space 
have rational values. If $f_1,\ldots,f_{2p}$ are locally transverse,
then the set of $(X,\si,\pi)$ giving doubly-periodic $\psi$ will
be countable and dense in ${\cal M}_p$. Having fixed $(X,\si,\pi)$ 
there are $p-2$ degrees of freedom to choose $\cal L$, all of which 
give doubly-periodic~$\psi$.

Next we do a similar parameter count for minimal Lagrangian tori. To 
do this we have to include the holomorphic involution $\rho:X\ra X$
as in \S\ref{is42}. This has the property that if $\pi(x)=\la$ then 
$\pi\circ\rho(x)=-\la$. Clearly, $\rho$ takes branch points to branch
points, so the set $\{\la_1,\ldots,\la_{2p}\}$ is invariant under
$\la\mapsto-\la$. It follows that $p$ is even, say $p=2d$, and we
can order the $\la_k$ such that $\la_{d+k}=-\la_k$ for~$k=1,\ldots,d$.

Thus, $\la_1,\ldots,\la_{4d}$ are determined by $\la_1,\ldots,\la_d$,
so there are $2d$ real parameters in $(X,\rho,\si,\pi)$. Suppose that
$d\ge 2$, so that $p\ge 4$. (If $d=1$ the solutions have an $\R$
symmetry group, and the parameter count is slightly different.)
The condition for the Tzitz\'eica solution to be doubly-periodic is 
$2d-4$ rationality conditions, and for $\psi$ to be doubly-periodic 
in $\CP^2$ is another 4 rationality conditions. As $L$ must be 
invariant under $\rho$ and $\si$, it has $d-2$ degrees of freedom.

Effectively, this means that the moduli space ${\cal M}_{2d}'$ of data
$(X,\rho,\si,\pi)$ has dimension $2d$, and to be doubly-periodic requires 
that $2d$ real functions $f_1,\ldots,f_{2d}$ on this moduli space have 
rational values. Having fixed $(X,\si,\pi)$ there are $d-2$ degrees of 
freedom to choose $\cal L$, all of which give doubly-periodic~$\psi$.

Here are our conclusions in brief:
\begin{itemize} 
\item Up to isometries of $\CP^2$, we expect the family of minimal 
tori in $\CP^2$ with spectral curve of genus $p\ge 2$ to depend on 
$2p$ rational numbers and $p-2$ real numbers.
\item Up to isometries of $\CP^2$, we expect the family of minimal 
Lagrangian tori in $\CP^2$ with spectral curve of genus $2d\ge 4$ 
to depend on $2d$ rational numbers and $d-2$ real numbers.
\end{itemize}
As the double-periodicity conditions for Lagrangian $\psi:\R^2\ra\CP^2$ 
are equivalent to those for its Legendrian lift $\phi:\R^2\ra{\cal S}^5$,
the second parameter count also gives the answer for the families of 
minimal Legendrian tori in $\cal S^5$ up to transformations in $\U(3)$,
and of special Lagrangian $T^2$-cones in $\C^3$ up to transformations 
in~$\SU(3)$.

The moral for special Lagrangian geometry is that one should expect
very large numbers of SL $T^2$-cones in $\C^3$, which can even exist
in continuous families up to isomorphisms. These provide many local
models for singularities of SL 3-folds in Calabi--Yau 3-folds.

\section{A family of special Lagrangian cones in $\C^3$}
\label{is5}

In \cite[\S 8]{Joyc1} and \cite[\S 6]{Joyc2} the author gave two 
constructions of countable families of special Lagrangian $T^2$-cones 
in $\C^3$, the first using $\U(1)$-invariance, and the second by 
evolving a 1-parameter family of quadric cones in Lagrangian planes.
Both constructions are related to work of other authors. 

In particular, the section on $\U(1)$-invariant SL cones in 
\cite{Joyc1} essentially repeats the work of Castro and Urbano 
\cite{CaUr} on $\U(1)$-invariant minimal tori in $\CP^2$, and 
was also discovered independently by Haskins \cite{Hask}. The 
`evolving quadrics' construction of \cite{Joyc2} generalizes 
examples of Lawlor \cite{Lawl}, and some of the examples it 
produces were also studied by Bryant \cite[\S 3.5]{Brya} from a 
different point of view. For more details, see~\cite{Joyc1,Joyc2}.

Motivated by these we shall now construct a more general family 
of special Lagrangian cones in $\C^3$ which includes those of
\cite{Joyc1,Joyc2} as special cases. These some from a family of 
explicit conformal harmonic maps $\phi:\R^2\ra{\cal S}^5$ with 
Legendrian image, which will be analyzed from the integrable
systems point of view in~\S\ref{is6}. 

\subsection{Constructing the family}
\label{is51}

Here is our main result.

\begin{thm} Let\/ $\be_1,\be_2,\be_3$ and\/ $\ga_1,\ga_2,\ga_3$ be real
numbers with not all\/ $\be_j$ and not all\/ $\ga_j$ zero, such that
\e
\be_1+\be_2+\be_3=0,\quad
\ga_1+\ga_2+\ga_3=0\quad\text{and}\quad
\be_1\ga_1+\be_2\ga_2+\be_3\ga_3=0.
\label{is5eq1}
\e
Suppose $y_1,y_2,y_3:\R\ra\C$ and\/ $v:\R\ra\R$ are functions of\/ $s$,
and\/ $z_1,z_2,z_3:\R\ra\C$ and\/ $w:\R\ra\R$ functions of\/ $t$, satisfying
\begin{align}
\frac{\d y_1}{\d s}&=\be_1\,\overline{y_2y_3\!}\,,&\quad
\frac{\d y_2}{\d s}&=\be_2\,\overline{y_3y_1\!}\,,&\quad
\frac{\d y_3}{\d s}&=\be_3\,\overline{y_1y_2\!}\,,
\label{is5eq2}\\
\frac{\d z_1}{\d t}&=\ga_1\,\overline{z_2z_3\!}\,,&\quad
\frac{\d z_2}{\d t}&=\ga_2\,\overline{z_3z_1\!}\,,&\quad
\frac{\d z_3}{\d t}&=\ga_3\,\overline{z_1z_2\!}\,,
\label{is5eq3}\\
\ms{y_1}&=\be_1v+1,&\quad
\ms{y_2}&=\be_2v+1,&\quad
\ms{y_3}&=\be_3v+1,
\label{is5eq4}\\
\ms{z_1}&=\ga_1w+1,&\quad
\ms{z_2}&=\ga_2w+1,&\quad
\ms{z_3}&=\ga_3w+1.
\label{is5eq5}
\end{align}
If\/ \eq{is5eq2}--\eq{is5eq3} hold for all\/ $s,t$ and\/
\eq{is5eq4}--\eq{is5eq5} hold for $s=t=0$, then \eq{is5eq4}--\eq{is5eq5}
hold for all\/ $s,t$, for some functions $v,w$. Define $\Phi:\R^3\ra\C^3$~by
\e
\Phi:(r,s,t)\mapsto{\ts\frac{1}{\sqrt{3}}}
\bigl(ry_1(s)z_1(t),ry_2(s)z_2(t),ry_3(s)z_3(t)\bigr).
\label{is5eq6}
\e
Define a subset\/ $N$ of\/ $\C^3$ by
\e
N=\bigl\{\Phi(r,s,t):r,s,t\in\R\bigr\}.
\label{is5eq7}
\e
Then $N$ is a special Lagrangian cone in~$\C^3$.
\label{is5thm1}
\end{thm}

\begin{proof} Suppose \eq{is5eq2} holds for all $s$, and \eq{is5eq4}
for $s=0$. From \eq{is5eq2} we deduce that $\frac{\d}{\d s}\bigl(
\ms{y_j}\bigr)=2\be_j\Re(y_1y_2y_3)$ for $j=1,2,3$. Comparing this
with \eq{is5eq4} shows that $v(s)$ should satisfy $\frac{\d v}{\d s}
=2\Re(y_1y_2y_3)$. Therefore, setting $v(s)=v(0)+2\int_0^s\Re
\bigl(y_1(u)y_2(u)y_3(u)\bigr)\d u$ shows that \eq{is5eq4} holds 
for all $s$. Similarly, if \eq{is5eq3} holds for all $t$ and 
\eq{is5eq5} holds for $t=0$, then it holds for all $t$. This 
proves the first part of the theorem. 

For the second part, we must show that $N$ is special Lagrangian 
wherever it is nonsingular, that is, wherever $\Phi$ is an immersion.
Now $\Phi$ is an immersion at $(r,s,t)$ when $\frac{\pd\Phi}{\pd r},
\frac{\pd\Phi}{\pd s},\frac{\pd\Phi}{\pd t}$ are linearly independent,
and then $T_{\Phi(r,s,t)}N=\ban{\frac{\pd\Phi}{\pd r},\frac{\pd\Phi}{\pd s},
\frac{\pd\Phi}{\pd t}}_{\R}$. Thus we must show that $T_{\Phi(r,s,t)}N$
is an SL 3-plane $\R^3$ in $\C^3$ for all $(r,s,t)$ for which
$\frac{\pd\Phi}{\pd r},\frac{\pd\Phi}{\pd s},\frac{\pd\Phi}{\pd t}$ are
linearly independent. By Proposition \ref{is2prop}, this holds if and only if
\ea
\om\Bigl(\frac{\pd\Phi}{\pd r},\frac{\pd\Phi}{\pd s}\Bigr)\equiv
\om\Bigl(\frac{\pd\Phi}{\pd r},\frac{\pd\Phi}{\pd t}\Bigr)\equiv
\om\Bigl(\frac{\pd\Phi}{\pd s},\frac{\pd\Phi}{\pd t}\Bigr)&\equiv 0
\label{is5eq8}\\
\text{and}\qquad \Im\Om\Bigl(\frac{\pd\Phi}{\pd r},
\frac{\pd\Phi}{\pd s},\frac{\pd\Phi}{\pd t}\Bigr)&\equiv 0.
\label{is5eq9}
\ea

Using equations \eq{is5eq2}, \eq{is5eq3} and \eq{is5eq6} we find that
\ea
\frac{\pd\Phi}{\pd r}&={\ts\frac{1}{\sqrt{3}}}
\bigl(y_1z_1,y_2z_2,y_3z_3\bigr),
\label{is5eq10}\\
\frac{\pd\Phi}{\pd s}&={\ts\frac{1}{\sqrt{3}}}
\bigl(r\be_1\,\overline{y_2y_3\!}\,z_1,
r\be_2\,\overline{y_3y_1\!}\,z_2,r\be_3\,\overline{y_1y_2\!}\,z_3\bigr),
\label{is5eq11}\\
\frac{\pd\Phi}{\pd t}&={\ts\frac{1}{\sqrt{3}}}
\bigl(r\ga_1y_1\,\overline{z_2z_3\!}\,,
r\ga_2y_2\,\overline{z_3z_1\!}\,,r\ga_3y_3\,\overline{z_1z_2\!}\,\bigr).
\label{is5eq12}
\ea
From \eq{is2eq} we deduce that $\om\bigl((a_1,a_2,a_3),(b_1,b_2,b_3)\bigr)
=\Im(a_1\bar b_1+a_2\bar b_2+a_3\bar b_3\bigr)$. Thus from \eq{is5eq10}
and \eq{is5eq11} we have
\begin{align*}
\om\bigl({\ts\frac{\pd\Phi}{\pd r},\frac{\pd\Phi}{\pd s}}\bigr)
&={\ts\frac{1}{3}}
r\Im(y_1y_2y_3)\bigl(\be_1\ms{z_1}+\be_2\ms{z_2}+\be_3\ms{z_3}\bigr)\\
&={\ts\frac{1}{3}}
r\Im(y_1y_2y_3)\bigl(\be_1(\ga_1w+1)+\be_2(\ga_2w+1)+\be_3(\ga_3w+1)\bigr)\\
&={\ts\frac{1}{3}}
r\Im(y_1y_2y_3)\bigl(\be_1+\be_2+\be_3+w(\be_1\ga_1+\be_2\ga_2+\be_3\ga_3)
\bigr)=0,
\end{align*}
using \eq{is5eq5} in the second line and \eq{is5eq1} in the third. This
proves the first equation of \eq{is5eq8}. The second follows in the same
way, and the third from
\begin{equation*}
\om\bigl({\ts\frac{\pd\Phi}{\pd s},\frac{\pd\Phi}{\pd t}}\bigr)=
{\ts\frac{1}{3}}r^2\Im(\,\overline{y_1y_2y_3\!}\,z_1z_2z_3)
(\be_1\ga_1+\be_2\ga_2+\be_3\ga_3)=0,
\end{equation*}
using \eq{is5eq1}, \eq{is5eq11} and \eq{is5eq12}.

To prove \eq{is5eq9}, observe that
\begin{align*}
\Om\Bigl(\frac{\pd\Phi}{\pd r}&,\frac{\pd\Phi}{\pd s},
\frac{\pd\Phi}{\pd t}\Bigr)=
\Bigl\vert\frac{\pd\Phi}{\pd r}\,\,\frac{\pd\Phi}{\pd s}\,\,
\frac{\pd\Phi}{\pd t}\Bigr\vert=
\frac{1}{3\sqrt{3}}
\left\vert\begin{matrix}
y_1z_1 & r\be_1\,\overline{y_2y_3\!}\,z_1 & r\ga_1y_1\,\overline{z_2z_3\!}\, \\
y_2z_2 & r\be_2\,\overline{y_3y_1\!}\,z_2 & r\ga_2y_2\,\overline{z_3z_1\!}\, \\
y_3z_3 & r\be_3\,\overline{y_1y_2\!}\,z_3 & r\ga_3y_3\,\overline{z_1z_2\!}\, 
\end{matrix}\right\vert \\
={\ts\frac{1}{3\sqrt{3}}}r^2\bigl(&\be_2\ms{y_3y_1}\ga_3\ms{z_1z_2}
+\be_3\ms{y_1y_2}\ga_1\ms{z_2z_3}+\be_1\ms{y_2y_3}\ga_2\ms{z_3z_1}\\
-&\be_3\ms{y_1y_2}\ga_2\ms{z_3z_1}-\be_1\ms{y_2y_3}
\ga_3\ms{z_1z_2}-\be_2\ms{y_3y_1}\ga_1\ms{z_2z_3}\bigr),
\end{align*}
where in the first line the terms $\md{\ldots}$ are determinants of
complex $3\t 3$ matrices, and $\frac{\pd\Phi}{\pd r},\frac{\pd\Phi}{\pd s},
\frac{\pd\Phi}{\pd t}$ are regarded are complex column matrices. As every
term in the second line is real, $\Om\bigl(\frac{\pd\Phi}{\pd r},
\frac{\pd\Phi}{\pd s},\frac{\pd\Phi}{\pd t}\bigr)$ is real. Thus
$\Im\Om\bigl(\frac{\pd\Phi}{\pd r},\frac{\pd\Phi}{\pd s},
\frac{\pd\Phi}{\pd t}\bigr)=0$, proving~\eq{is5eq9}.
\end{proof}

\subsection{Explicit solution of the o.d.e.s \eq{is5eq2} and \eq{is5eq3}} 
\label{is52}

As in \cite[\S 8]{Joyc1} and \cite[\S 6]{Joyc2}, we can simplify the 
solutions of \eq{is5eq2} and \eq{is5eq3}. To do this, note first that 
\eq{is5eq2} implies that
\begin{equation*}
{\ts\frac{\d}{\d s}}(y_1y_2y_3)=\be_1\ms{y_2y_3}+\be_2\ms{y_3y_1}
+\be_3\ms{y_1y_2}.
\end{equation*}
As the right hand side is real, we have $\Im(y_1y_2y_3)\equiv B$ 
for some~$B\in\R$.

Now $\ms{y_1y_2y_3}=(\be_1v+1)(\be_2v+1)(\be_3v+1)$, so this gives
\begin{equation*}
\bms{\Re(y_1y_2y_3)}=(\be_1v+1)(\be_2v+1)(\be_3v+1)-B^2.
\end{equation*}
However, $\frac{\d v}{\d s}=2\Re(y_1y_2y_3)$, as in the proof 
of Theorem \ref{is5thm1}, and so $v$ satisfies the o.d.e.
\begin{equation*}
\Bigl(\frac{\d v}{\d s}\Bigr)^2=4\bigl((\be_1v+1)(\be_2v+1)(\be_3v+1)
-B^2\bigr).
\end{equation*}
Since $\ms{y_j}=\be_jv+1$ by \eq{is5eq4} we may write
$y_j(s)={\rm e}^{i\de_j(s)}\sqrt{\be_jv(s)\!+\!1}$ for $j=1,2,3$,
for real functions $\de_1,\de_2,\de_3$. In this way we prove:

\begin{prop} In the situation of Theorem \ref{is5thm1} the 
functions $y_1,y_2,y_3$ may be written $y_j(s)={\rm e}^{i\de_j(s)}
\sqrt{\be_jv(s)\!+\!1}$, for $v,\de_1,\de_2,\de_3:\R\ra\R$. Define 
$Q(v)=(\be_1v+1)(\be_2v+1)(\be_3v+1)$ and\/ $\de=\de_1+\de_2+\de_3$.
Then \eq{is5eq4} holds automatically, and\/ \eq{is5eq2} is equivalent to
\e
\Bigl(\frac{\d v}{\d s}\Bigr)^2=4\bigl(Q(v)-B^2\bigr)
\quad\text{and}\quad 
\frac{\d\de_j}{\d s}=-\,\frac{\be_jB}{\be_jv+1}
\label{is5eq13}
\e
for $j=1,2,3$, where $\Im(y_1y_2y_3)\equiv Q(v)^{1/2}\sin\de\equiv B$ 
for some~$B\in[-1,1]$. 
\label{is5prop1}
\end{prop}

Here is the corresponding result for $z_1,z_2,z_3$.

\begin{prop} In the situation of Theorem \ref{is5thm1} the 
functions $z_1,z_2,z_3$ may be written $z_j(t)={\rm e}^{i\ep_j(t)}
\sqrt{\ga_jw(t)\!+\!1}$, for $w,\ep_1,\ep_2,\ep_3:\R\ra\R$. Define 
$R(w)=(\ga_1w+1)(\ga_2w+1)(\ga_3w+1)$ and\/ $\ep=\ep_1+\ep_2+\ep_3$.
Then \eq{is5eq5} holds automatically, and\/ \eq{is5eq3} is equivalent to
\e
\Bigl(\frac{\d w}{\d t}\Bigr)^2=4\bigl(R(w)-C^2\bigr)
\quad\text{and}\quad 
\frac{\d\ep_j}{\d t}=-\,\frac{\ga_jC}{\ga_jw+1}
\label{is5eq14}
\e
for $j=1,2,3$, where $\Im(z_1z_2z_3)\equiv R(w)^{1/2}\sin\ep\equiv C$ 
for some~$C\in[-1,1]$. 
\label{is5prop2}
\end{prop}

As in \cite[\S 8.2]{Joyc1}, the o.d.e.s for $v$ and $w$ in 
\eq{is5eq13} and \eq{is5eq14} can be solved entirely explicitly 
in terms of the {\it Jacobi elliptic functions}. Then $\de_j$ 
and $\ep_j$ can also be given explicitly, in terms of integrals 
involving the Jacobi elliptic functions, and so the solutions
$y_j,z_j$ of \eq{is5eq2} and \eq{is5eq3} are known explicitly
in terms of elliptic functions. We shall not give these solutions 
here.

\subsection{Conformal harmonic maps to ${\cal S}^5$ and $\CP^2$}
\label{is53}

As in \S\ref{is41}, a special Lagrangian cone in $\C^3$ induces
conformal harmonic maps $\phi:S\ra{\cal S}^5$ and $\psi:S\ra\CP^2$
from a Riemann surface $S$. We shall now write these out explicitly 
for the SL cones of Theorem \ref{is5thm1}. For convenience, we begin 
by choosing a normalization for the constants $\be_1,\be_2,\be_3$
and~$\ga_1,\ga_2,\ga_3$.

Regarding $\bs\be=(\be_1,\be_2,\be_3)$ and $\bs\ga=(\ga_1,\ga_2,\ga_3)$ 
as vectors in $\R^3$, the conditions on $\bs\be,\bs\ga$ in Theorem 
\ref{is5thm1} are that $\bs\be$ and $\bs\ga$ should be nonzero, and 
that $\bs\be,\bs\ga$ and $(1,1,1)$ should be orthogonal. However, 
multiplying $\bs\be$ or $\bs\ga$ by a nonzero constant has no effect on 
the set of special Lagrangian cones constructed in Theorem~\ref{is5thm1}. 

To see this, let $\be_j,\ga_j,y_j,z_j,v$ and $w$ satisfy the conditions 
of the theorem, let $\si,\tau\in\R$ be nonzero, and define
\e
\begin{aligned}
\be_j'&=\si\be_j,&\quad
y_j'(s)&=y_j(\si s)&\quad\text{for $j=1,2,3$, and}\quad
v'(s)&=\si^{-1}v(\si s),\\
\ga_j'&=\tau\ga_j,&\quad
z_j'(t)&=z_j(\tau t)&\quad\text{for $j=1,2,3$, and}\quad
w'(t)&=\tau^{-1}w(\tau t).
\end{aligned}
\label{is5eq15}
\e
Then it is easy to show that $\be_j',\ga_j',y_j',z_j',v'$ and $w'$
also satisfy the conditions of the theorem, yielding $\Phi':\R^3\ra\C^3$
with $\Phi'(r,s,t)=\Phi(r,\si s,\tau t)$, so that the images of $\Phi'$
and $\Phi$ are the same special Lagrangian cone. 

Therefore, we are free to rescale $\bs\be$ and $\bs\ga$ without 
changing the resulting set of SL cones. Fix $\md{\bs\be}=\md{\bs\ga}=1$. 
Then $\bs\be$ and $\bs\ga$ lie on the unit circle in the plane 
$x_1+x_2+x_3=0$ in $\R^3$ and are orthogonal, so we may write
\begin{align*}
\bs\be&=\cos\th\cdot{\ts\frac{1}{\sqrt{2}}}(1,-1,0)+
\sin\th\cdot{\ts\frac{1}{\sqrt{6}}}(-1,-1,2) \\
\text{and}\quad
\bs\ga&=\cos\th\cdot{\ts\frac{1}{\sqrt{6}}}(-1,-1,2) 
-\sin\th\cdot{\ts\frac{1}{\sqrt{2}}}(1,-1,0)
\end{align*}
for some $\th\in[0,2\pi)$. (Here $\bs\be$ determines $\bs\ga$ up to sign,
which we have chosen arbitrarily.) We shall show that with these
choices, the map $(s,t)\mapsto\Phi(1,s,t)$ is conformal.

\begin{thm} Fix $\th\in[0,2\pi)$, and define
\begin{align}
\be_1&\!=\!{\ts\frac{1}{\sqrt{2}}}\cos\th\!-\!
{\ts\frac{1}{\sqrt{6}}}\sin\th,&\;
\be_2&\!=\!-{\ts\frac{1}{\sqrt{2}}}\cos\th\!-\!
{\ts\frac{1}{\sqrt{6}}}\sin\th,&\;
\be_3&\!=\!{\ts\frac{2}{\sqrt{6}}}\sin\th,
\label{is5eq16}\\
\ga_1&\!=\!-{\ts\frac{1}{\sqrt{6}}}\cos\th\!-\!
{\ts\frac{1}{\sqrt{2}}}\sin\th,&\;
\ga_2&\!=\!-{\ts\frac{1}{\sqrt{6}}}\cos\th\!+\!
{\ts\frac{1}{\sqrt{2}}}\sin\th,&\;
\ga_3&\!=\!{\ts\frac{2}{\sqrt{6}}}\cos\th.
\label{is5eq17}
\end{align}
In the situation of Theorem \ref{is5thm1}, with these values of\/
$\be_j$ and\/ $\ga_j$, we have $\ms{\Phi(r,s,t)}=r^2$ and
$\frac{\pd\Phi}{\pd r}$, $\frac{\pd\Phi}{\pd s}$ and\/
$\frac{\pd\Phi}{\pd t}$ are orthogonal with\/ 
$\ms{\frac{\pd\Phi}{\pd r}}=1$ and
\ea
\bms{\ts\frac{\pd\Phi}{\pd s}}&=\bms{\ts\frac{\pd\Phi}{\pd t}}
=2r^2\bigl(a+bv(s)+cw(t)\bigr),
\label{is5eq18}\\
\text{where}\quad
a&={\ts\frac{1}{6}}(\be_1^2+\be_2^2+\be_3^2)=
{\ts\frac{1}{6}}(\ga_1^2+\ga_2^2+\ga_3^2)={\ts\frac{1}{6}},
\label{is5eq19}\\
b&=-{\ts\frac{1}{6}}(\be_1^3\!+\!\be_2^3\!+\!\be_3^3)=
{\ts\frac{1}{6}}(\be_1\ga_1^2\!+\!\be_2\ga_2^2\!+\!\be_3\ga_3^2)
=-\ha\be_1\be_2\be_3,
\label{is5eq20}\\
\text{and}\quad 
c&={\ts\frac{1}{6}}(\be_1^2\ga_1\!+\!\be_2^2\ga_2\!+\!\be_3^2\ga_3)
=-{\ts\frac{1}{6}}(\ga_1^3\!+\!\ga_2^3\!+\!\ga_3^3)
=-\ha\ga_1\ga_2\ga_3.
\label{is5eq21}
\ea
The maps $\phi:\R^2\ra{\cal S}^5$ and\/ $\psi:\R^2\ra\CP^2$
defined by $\phi:(s,t)\mapsto\Phi(1,s,t)$ and\/ $\psi:(s,t)
\mapsto[\Phi(1,s,t)]$ are both conformal harmonic maps.
\label{is5thm2}
\end{thm}

\begin{proof} For the first part, by \eq{is5eq1} and
\eq{is5eq4}--\eq{is5eq6} we have
\begin{align*}
\bms{\Phi(r,s,t)}&={\ts\frac{1}{3}}r^2\bigl(\ms{y_1}\ms{z_1}+
\ms{y_2}\ms{z_2}+\ms{y_3}\ms{z_3}\bigr)\\
&={\ts\frac{1}{3}}r^2\bigl((\be_1v\!+\!1)(\ga_1w\!+\!1)\!+\!
(\be_2v\!+\!1)(\ga_2w\!+\!1)\!+\!(\be_3v\!+\!1)(\ga_3w\!+\!1)\bigr)\\
&={\ts\frac{1}{3}}r^2\bigl(3\!+\!(\be_1\!+\!\be_2\!+\!\be_3)v\!+
\!(\ga_1\!+\!\ga_2\!+\!\ga_3)w\!+\!(\be_1\ga_1\!+\!\be_2\ga_2\!+
\!\be_3\ga_3)vw\bigr)\\
&=r^2.
\end{align*}
The equation $\ms{\frac{\pd\Phi}{\pd r}}=1$ follows in the same way.

To prove $\frac{\pd\Phi}{\pd r},\frac{\pd\Phi}{\pd s},\frac{\pd\Phi}{\pd
t}$ are orthogonal we use \eq{is5eq10}--\eq{is5eq12} and the formula
\hfil\break
$g\bigl((a_1,a_2,a_3),(b_1,b_2,b_3)\bigr)=\Re(a_1\bar b_1+a_2\bar b_2
+a_3\bar b_3\bigr)$. Thus we have
\begin{align*}
g\bigl({\ts\frac{\pd\Phi}{\pd r},\frac{\pd\Phi}{\pd s}}\bigr)
&={\ts\frac{1}{3}}r\Re(y_1y_2y_3)\bigl(\be_1\ms{z_1}+
\be_2\ms{z_2}+\be_3\ms{z_3}\bigr)\\
&={\ts\frac{1}{3}}r\Re(y_1y_2y_3)\bigl(\be_1(\ga_1w+1)
+\be_2(\ga_2w+1)+\be_3(\ga_3w+1)\bigr)\\
&={\ts\frac{1}{3}}r\Re(y_1y_2y_3)\bigl(\be_1+\be_2+\be_3
+w(\be_1\ga_1+\be_2\ga_2+\be_3\ga_3)\bigr)=0,
\end{align*}
using \eq{is5eq5} in the second line and \eq{is5eq1} in the third. In
the same way we show that $g\bigl(\frac{\pd\Phi}{\pd r},\frac{\pd\Phi}{\pd
t}\bigr)=g\bigl(\frac{\pd\Phi}{\pd s},\frac{\pd\Phi}{\pd t}\bigr)=0$, and
so $\frac{\pd\Phi}{\pd r},\frac{\pd\Phi}{\pd s},\frac{\pd\Phi}{\pd t}$ are
orthogonal.

Using equations \eq{is5eq4}, \eq{is5eq5} and \eq{is5eq11} we obtain
\begin{align*}
\bms{\ts\frac{\pd\Phi}{\pd s}}&={\ts\frac{1}{3}}r^2\bigl[
\be_1^2\ms{y_2}\ms{y_3}\ms{z_1}+
\be_2^2\ms{y_3}\ms{y_1}\ms{z_2}+
\be_3^2\ms{y_1}\ms{y_2}\ms{z_3}\bigr]\\
&={\ts\frac{1}{3}}r^2\bigl[\be_1^2(\be_2v\!+\!1)(\be_3v\!+\!1)
(\ga_1w\!+\!1)\!+\!\be_2^2(\be_3v\!+\!1)(\be_1v\!+\!1)(\ga_2w\!+\!1)\\
&\qquad\quad+\!\be_3^2(\be_1v\!+\!1)(\be_2v\!+\!1)(\ga_3w\!+\!1)\bigr]\\
&={\ts\frac{1}{3}}r^2\bigl[(\be_1^2\!+\!\be_2^2\!+\!\be_3^2)\!+\!
v\bigl(\be_1^2(\be_2\!+\!\be_3)\!+\!\be_2^2(\be_3\!+\!\be_1)\!+\!
\be_3^2(\be_1\!+\!\be_2)\bigr)\\
+\!w(&\be_1^2\ga_1\!+\!\be_2^2\ga_2\!+\!\be_3^2\ga_3)
\!+\!vw\bigl(\be_1^2(\be_2\!+\!\be_3)\ga_1\!+\!\be_2^2(\be_3\!+\!\be_1)
\ga_2\!+\!\be_3^2(\be_1\!+\!\be_2)\ga_3\bigr)\\
&\qquad\quad+\!v^2\be_1\be_2\be_3(\be_1\!+\!\be_2\!+\!\be_3)
\!+\!v^2w\be_1\be_2\be_3(\be_1\ga_1\!+\!\be_2\ga_2\!+\!\be_3\ga_3)\bigr].
\end{align*}

By \eq{is5eq1}, the terms in $v^2$ and $v^2w$ vanish. Also, using 
\eq{is5eq16} and \eq{is5eq17} we have
\begin{align*}
\be_1^2(\be_2&+\be_3)\ga_1+\be_2^2(\be_3+\be_1)
\ga_2+\be_3^2(\be_1+\be_2)\ga_3\\
&={\ts\frac{1}{\sqrt{3}}}\bigl(\be_1^2(\be_2^2-\be_3^2)+
\be_2^2(\be_3^2-\be_1^2)+\be_3^2(\be_1^2-\be_2^2)\bigr)=0,
\end{align*}
so the term in $vw$ vanishes. Thus, replacing $(\be_2+\be_3)$ by 
$-\be_1$, etc., we get
\begin{equation*}
\bms{\ts\frac{\pd\Phi}{\pd s}}={\ts\frac{1}{3}}
r^2\bigl[(\be_1^2+\be_2^2+\be_3^2)-
v(\be_1^3+\be_2^3+\be_3^3)
+w(\be_1^2\ga_1+\be_2^2\ga_2+\be_3^2\ga_3)\bigr].
\end{equation*}
In the same way, we find that
\begin{equation*}
\bms{\ts\frac{\pd\Phi}{\pd t}}
={\ts\frac{1}{3}}r^2\bigl[(\ga_1^2+\ga_2^2+\ga_3^2)
+v(\be_1\ga_1^2+\be_2\ga_2^2+\be_3\ga_3^2)
-w(\ga_1^3+\ga_2^3+\ga_3^3)\bigr].
\end{equation*}

Now using \eq{is5eq16} and \eq{is5eq17} one can show that
\begin{gather*}
\be_1^2\!+\!\be_2^2\!+\!\be_3^2=\ga_1^2\!+\!\ga_2^2\!+\!\ga_3^2,\quad
-(\be_1^3\!+\!\be_2^3\!+\!\be_3^3)=
\be_1\ga_1^2\!+\!\be_2\ga_2^2\!+\!\be_3\ga_3^2=-3\be_1\be_2\be_3\\
\text{and}\quad \be_1^2\ga_1\!+\!\be_2^2\ga_2\!+\!\be_3^2\ga_3=
-(\ga_1^3\!+\!\ga_2^3\!+\!\ga_3^3)=-3\ga_1\ga_2\ga_3.
\end{gather*}
The last five equations prove \eq{is5eq18}--\eq{is5eq21}, as we 
want. Finally, it follows from what we have proved so far that 
$\phi:\R^2\ra{\cal S}^5$ is a conformal map, and as its image is 
minimal, it is also harmonic. As $\phi$ has Legendrian image, 
$\psi$ is also conformal and harmonic, in the usual way.
\end{proof}

As from \S\ref{is52} the functions $y_j,z_j$ defining $\Phi$ are 
known explicitly in terms of integrals involving the Jacobi elliptic 
functions, we have constructed families of {\it explicit conformal 
harmonic maps} $\phi:\R^2\ra{\cal S}^5$ and~$\psi:\R^2\ra\CP^2$.

\subsection{Interesting special cases, and double periodicity}
\label{is54}

We now consider some special cases in which the $y_j$ or $z_j$
assume a simple form, and so explain how to recover the constructions
of \cite[\S 8]{Joyc1} and \cite[\S 6]{Joyc2} from the more general
construction above.

\begin{itemize}
\item[(a)] Let $\ka_1,\ka_2,\ka_3\in\R$ with $\ka_1+\ka_2+\ka_3=-\pi/2$,
and define $y_j={\rm e}^{i(\be_js+\ka_j)}$ for $j=1,2,3$. Then it is easy
to see that $y_1,y_2,y_3$ satisfy \eq{is5eq2} and \eq{is5eq4}, with
$v\equiv 0$ and $B=-1$. The corresponding special Lagrangian cones in 
Theorem \ref{is5thm1} are invariant under the group action
\begin{equation*}
(z_1,z_2,z_3)\mapsto({\rm e}^{i\be_1s}z_1,{\rm e}^{i\be_2s}z_2,
{\rm e}^{i\be_3s}z_3)
\end{equation*}
for $s\in\R$, which is a $\U(1)$ subgroup of $\SU(3)$ if 
$\be_1,\be_2,\be_3$ are relatively rational, and an $\R$ 
subgroup otherwise. In this case, Theorem \ref{is5thm1} 
reduces to the construction of $\U(1)$-invariant SL cones 
in $\C^3$ given in~\cite[\S 8]{Joyc1}. 

In the same way, we can take $z_j={\rm e}^{i(\ga_jt+\ka_j)}$ for 
$j=1,2,3$, with $w\equiv 0$ and $C=-1$, and two similar cases with
$B=1$ and $C=1$, all of which give $\U(1)$-invariant or $\R$-invariant
SL cones in $\C^3$ coming from the construction of \cite[\S 8]{Joyc1}.
(See also Castro and Urbano \cite{CaUr}, and Haskins~\cite{Hask}.)

\item[(b)] Take $B=0$ in Proposition \ref{is5prop1}. Then \eq{is5eq13}
shows that the phases $\de_1,\de_2,\de_3$ are constant, so we may as
well fix them to be 0 or $\pi$, and take $y_1,y_2,y_3$ to be {\it real}.
As in \cite[\S 6.1]{Joyc2} the $y_j$ are given by simple formulae
involving Jacobi elliptic functions (rather than integrals of Jacobi
elliptic functions). 

As the point $(y_1,y_2,y_3)$ moves in $\R^3$ it sweeps out one of the
two connected components of the curve
\begin{equation*}
\bigl\{(x_1,x_2,x_3)\in\R^3:x_1^2+x_2^2+x_3^2=3,\quad
\ga_1x_1^2+\ga_2x_2^2+\ga_3x_3^2=0\bigr\}.
\end{equation*}
From this it follows that for fixed $t$, as $r,s$ vary $\Phi(r,s,t)$ 
sweeps out a quadric cone in a Lagrangian $\R^3$ in $\C^3$. So the 
special Lagrangian cone $N$ of \eq{is5eq7} is the total space of
a 1-parameter family of such quadrics, and we recover the `evolving 
quadrics' construction of~\cite{Joyc2}. 

In the same way, if $C=0$ in Proposition \ref{is5prop2} a similar 
thing happens, with $s$ and $t$ exchanged.

\item[(c)] Set $\th=0$ in Theorem \ref{is5thm2}. Then $\be_3=0$,
so $y_3$ is constant with $\md{y_3}=1$ by \eq{is5eq2} and \eq{is5eq4},
and $y_1,y_2$ are linear combinations of ${\rm e}^{\pm is/\sqrt{2}}$. 
Also $\ga_1=\ga_2$, so $z_2\equiv{\rm e}^{i\ka}z_1$ for some~$\ka\in\R$.

The corresponding SL cones in $\C^3$ turn out to be invariant
under a $\U(1)$ subgroup of $\SU(3)$ which fixes the third
coordinate in $\C^3$, corresponding to translation in the $s$
variable. Thus, after a linear coordinate change in $\C^3$,
this case reduces to a special case of the $\U(1)$-invariant
cones in part (a), but with a different parametrization.

In the same way, for each of the five other values of $\th\in[0,2\pi)$
for which one of $\be_2,\be_3,\ga_1,\ga_2$ and $\ga_3$ is zero, a
similar thing happens.
\end{itemize}

Next we consider when the maps $\phi:\R^2\ra{\cal S}^5$ and $\psi:
\R^2\ra\CP^2$ of Theorem \ref{is5thm2} are {\it doubly-periodic} in 
$\R^2$. Then $\phi$ and $\psi$ push down to conformal harmonic maps
$T^2\ra{\cal S}^5$ and $T^2\ra\CP^2$ whose images are {\it minimal 
tori} in ${\cal S}^5$ and $\CP^2$, and the special Lagrangian cone
$N$ of Theorem \ref{is5thm1} is a cone on $T^2$. We suppose for
simplicity that $\be_j,\ga_j$ are normalized as in 
equations~\eq{is5eq16}--\eq{is5eq17}.

It turns out that in cases (a)--(c) above the double-periodicity
conditions are soluble:
\begin{itemize}
\item[(a)] In case (a), suppose $\be_1,\be_2,\be_3$ are 
{\it relatively rational}. This happens for a countable dense
set of $\th\in[0,2\pi)$. Then $\be_j=n_j/S$ for $S>0$ and 
$n_1,n_2,n_3$ coprime integers. It follows that $y_j(s+S)=y_j(s)$ 
for $j=1,2,3$ and $s\in\R$, so that the $y_j$ are periodic in~$s$.

For double periodicity in $s,t$, the $z_j$ have only to be 
periodic up to multiplication by ${\rm e}^{i\be_js}$ for some 
$s\in\R$. Now $\th$ and $B$ are already fixed, but we are free 
to vary the constant $C$ in Proposition \ref{is5prop2}. It is 
shown in \cite[Th.~8.5]{Joyc1} that double periodicity holds 
for a countable dense subset of~$C\in[-1,1]$.
\item[(b)] In case (b) with $B=0$, $y_1,y_2,y_3$ are automatically 
periodic in $s$. We then need to vary the remaining data $\th,C$
to make $z_1,z_2,z_3$ periodic in~$t$.

Now $w$ is always periodic in $t$, with period $T$, say, and
the $z_j$ transform as $z_j(t+T)={\rm e}^{i\ze_j}z_j(t)$, where 
$\ze_1,\ze_2,\ze_3\in\R$ with $\ze_1+\ze_2+\ze_3=0$. If 
$\ze_j\in\pi\Q$ for $j=1,2,3$ then $n\ze_j\in 2\pi\Z$ for some 
positive integer $n$, and then $z_1,z_2,z_3$ are periodic with 
period $nT$. So, for double periodicity we need 2 functions of 
$\th$ and $C$ to be rational. In \cite[Th.s 5.9, 6.3 \& 6.4]{Joyc2}
it is shown that the $z_j$ are periodic for a countable dense
set of values of~$(\th,C)$.

\item[(c)] In case (c), $y_1,y_2,y_3$ are automatically periodic with 
period $2\sqrt{2}\pi$. Also, as $z_2\equiv{\rm e}^{i\ka}z_1$, the 
periodicity conditions for $z_1,z_2,z_3$ reduce to one rationality
condition, rather than two. As in case (a), $z_1,z_2,z_3$ are 
periodic in $t$ for a countable dense subset of~$C\in[-1,1]$.
\end{itemize}

What about double periodicity conditions in the general case?
If $\md{B}=1$ then $v$ is constant and we are in case (a) above,
so suppose $\md{B}<1$, and similarly $\md{C}<1$. Then $v,w$ are 
automatically nonconstant and periodic in $s,t$, with periods 
$S,T$ say, and the $y_j$ and $z_j$ transform as
\begin{equation*}
y_j(s+S)={\rm e}^{i\eta_j}y_j(s)\quad\text{and}\quad
z_j(t+T)={\rm e}^{i\ze_j}z_j(t)
\end{equation*}
for some constants $\eta_j,\zeta_j\in\R$ with $\eta_1+\eta_2+\eta_3
=\ze_1+\ze_2+\ze_3=0$. The conditions for the $y_j$ and $z_j$ to be
periodic in $s$ and $t$ are that $\eta_j\in\pi\Q$ and $\ze_j\in\pi\Q$
for $j=1,2,3$ respectively.

Thus, for $\phi$ and $\psi$ to be doubly-periodic we need 
the four functions $\eta_1/\pi$, $\eta_2/\pi$, $\ze_1/\pi$, 
$\ze_2/\pi$ of the three variables $\th,B,C$ to be rational. 
This is an {\it overdetermined}\/ problem, so it seems likely 
that in the general case, the double periodicity conditions 
will have few solutions, or none. Other than parts (a)--(c)
abive, the author knows of no cases in which $\phi,\psi$ are
doubly-periodic.

We can use \eq{is5eq18} to give a formula for the area of the minimal 
tori in ${\cal S}^5$ or $\CP^2$ arising from the construction above.

\begin{prop} Suppose that the map $\phi:\R^2\ra{\cal S}^5$ defined
in Theorem \ref{is5thm2} is doubly-periodic in $(s,t)$, with image 
$\Si$, so that\/ $\Si$ is a minimal torus in ${\cal S}^5$. Let\/ 
$S,T$ be the periods of\/ $v,w$ in $s$ and\/ $t$, as above, and 
let the period lattice of\/ $(s,t)\mapsto\Phi(1,s,t)$ in $\R^2$ 
be generated by $(a_{11}S,a_{12}T)$ and\/ $(a_{21}S,a_{22}T)$
for integers $a_{ij}$. Let\/ $N=\md{a_{11}a_{22}-a_{12}a_{21}}$. 
Then the area of\/ $\Si$ is
\e
{\mathop{\rm Area}}(\Si)=2N\Bigl(aST+bT\int_0^Sv(s)\d s
+cS\int_0^Tw(t)\d t\Bigr).
\label{is5eq23}
\e
\label{is5prop3}
\end{prop}

\begin{proof} As $\frac{\pd\Phi}{\pd s}$ and $\frac{\pd\Phi}{\pd t}$
are orthogonal, \eq{is5eq18} implies that the area form on $\Si$ is
$2\bigl(a+bv(s)+cw(t)\bigr)\d s\w\d t$. Also, as the period lattice
is generated by $(a_{11}S,a_{12}T)$ and $(a_{21}S,a_{22}T)$, we can
divide $\Si$ into $N=\md{a_{11}a_{22}-a_{12}a_{21}}$ copies of the
basic rectangle $[0,S]\t[0,T]$, each of which has area
$\int_0^T\int_0^S2\bigl(a+bv(s)+cw(t)\bigr)\d s\d t$. Equation
\eq{is5eq23} follows immediately.
\end{proof}

Observe that $v,w$ can be written explicitly using Jacobi elliptic
functions as in \S\ref{is52}, and so \eq{is5eq23} could easily be 
evaluated numerically in examples using a computer. This may be 
valuable in studying singularities of special Lagrangian 3-folds, 
since the area of $\Si$ is a crude measure of how nongeneric 
singularities modelled on the cone on $\Si$ are in the family 
of all special Lagrangian 3-folds. Also, note that the area 
of $\Si$ in ${\cal S}^5$ is the same as the area of its image 
in $\CP^2$, as the two are isometric.

\subsection{Comparison with constant mean curvature tori in $\R^3$}
\label{is55}

There is a strong analogy between the minimal Lagrangian tori 
in $\CP^2$ constructed above, and the examples of {\it constant 
mean curvature (CMC) tori} in $\R^3$ constructed by Wente \cite{Went} 
and Abresch \cite{Abre}, known as {\it Wente tori}. Wente proved 
\cite{Went} using analysis that there exist immersed CMC tori in $\R^3$, 
and so provided the first counterexamples to a conjecture of Hopf that 
the only compact surfaces in $\R^3$ with constant mean curvature are 
round spheres. 

Motivated by Wente's construction, Abresch \cite{Abre} gave explicit
formulae for the Wente tori in terms of elliptic integrals. Abresch's 
solutions are very similar in structure to those above. In particular, 
they have a `separated variable' form, being given in terms of 
single-variable functions $f(s),g(t)$ rather than two-variable 
functions, and $f$ and $g$ may be written explicitly using Jacobi 
elliptic functions.

We can also exploit the analogy in another way. Kapouleas \cite{Kapo} 
used analytic methods to construct examples of compact CMC surfaces 
$\Si$ in $\R^3$ for any genus $g\ge 3$. It seems very likely that one 
could use Kapouleas' method to construct examples of higher genus 
(immersed) minimal Lagrangian surfaces in $\CP^2$, and minimal 
Legendrian surfaces in~${\cal S}^5$. 

Kapouleas makes his examples by gluing together long segments of 
{\it Delaunay surfaces}, which are $\SO(2)$-invariant CMC surfaces 
resembling a string of round 2-spheres joined by narrow, catenoid-like 
`necks'. The appropriate analogues of Delaunay surfaces in our problem 
are Legendrian surfaces in ${\cal S}^5$ invariant under the $\U(1)$-action 
$(z_1,z_2,z_3)\mapsto({\rm e}^{is}z_1,{\rm e}^{-is}z_2,z_3)$, for~$s\in\R$.

In the notation of \S\ref{is51}--\S\ref{is53}, these have $\th=0$ 
and $B=-1$. When the remaining parameter $C\in[-1,1]$ is nonzero 
and small, the corresponding minimal Legendrian surfaces resemble 
chains of round Legendrian ${\cal S}^2$'s in ${\cal S}^5$ joined 
by small necks.

\section{Interpretation using integrable systems}
\label{is6}

In Theorem \ref{is5thm2} we constructed families of conformal harmonic 
maps $\phi:\R^2\ra{\cal S}^5$ and $\psi:\R^2\ra\CP^2$. We shall now 
analyze these in the integrable systems framework described in 
\S\ref{is3} and \S\ref{is4}. We will show that they are generically 
superconformal, and explicitly determine their harmonic sequences, 
Toda and Tzitz\'eica solutions, loops of flat connections, polynomial 
Killing fields, and spectral curves. This goes some way towards 
redressing the `dearth of examples' of superconformal harmonic tori 
referred to by Bolton and Woodward \cite[p.~76]{FoWo}. We shall use 
the notation of \S\ref{is51}--\S\ref{is53} throughout.

\subsection{The harmonic sequence of $\psi$}
\label{is61}

In the situation of \S\ref{is31}, take $U$ to be $\R^2$ with complex 
coordinate $z=s+it$. Then $\frac{\pd}{\pd z}=\ha\frac{\pd}{\pd s}-
\frac{i}{2}\frac{\pd}{\pd t}$ and $\frac{\pd}{\pd\bar z}=\ha
\frac{\pd}{\pd s}+\frac{i}{2}\frac{\pd}{\pd t}$. Thus by \eq{is5eq2}, 
\eq{is5eq3}, \eq{is5eq6} and the definition $\phi(s,t)=\Phi(1,s,t)$
we have
\begin{align*}
{\ts\frac{\pd\phi}{\pd z}}&=
{\ts\frac{1}{2\sqrt{3}}}
(\be_1\overline{y_2y_3\!}\,z_1\!-\!i\ga_1y_1\overline{z_2z_3\!}\,,
\be_2\overline{y_3y_1\!}\,z_2\!-\!i\ga_2y_2\overline{z_3z_1\!}\,,
\be_3\overline{y_1y_2\!}\,z_3\!-\!i\ga_3y_3\overline{z_1z_2\!}\,),\\
{\ts\frac{\pd\phi}{\pd\bar z}}&=
{\ts\frac{1}{2\sqrt{3}}}
(\be_1\overline{y_2y_3\!}\,z_1\!+\!i\ga_1y_1\overline{z_2z_3\!}\,,
\be_2\overline{y_3y_1\!}\,z_2\!+\!i\ga_2y_2\overline{z_3z_1\!}\,,
\be_3\overline{y_1y_2\!}\,z_3\!+\!i\ga_3y_3\overline{z_1z_2\!}\,).
\end{align*}
Calculation using \eq{is5eq4} and \eq{is5eq5} shows that
$\ban{\frac{\pd\phi}{\pd z},\phi}=\ban{\frac{\pd\phi}{\pd\bar z},\phi}=0$.
Also, using \eq{is5eq18} we find that~$\bms{\ts\frac{\pd\phi}{\pd z}}=
\bms{\ts\frac{\pd\phi}{\pd\bar z}}=a+bv(s)+cw(t)$.

As $\ban{\frac{\pd\phi}{\pd\bar z},\phi}=0$, by definition $\phi$ is
a holomorphic section of the holomorphic line bundle $L_0$ over $\C$
associated to $\psi_0=\psi:\C\ra\CP^2$. Therefore, from \S\ref{is31},
there exists a unique sequence of maps $\phi_k:\C\ra\C^3$ with
$\phi_0=\phi$, which satisfy \eq{is3eq1}, and the harmonic sequence
$(\psi_k)$ of $\psi$ is given by~$\psi_k=[\phi_k]$.

From \eq{is3eq1} we see that $\phi_{-1}=-\ms{\phi_0}
\md{\frac{\pd\phi_0}{\pd\bar z}}^{-2}\frac{\pd\phi_0}{\pd\bar z}$
and $\phi_1=\frac{\pd\phi_0}{\pd z}$, since $\md{\phi_0}\equiv 1$.
Thus the equations above give
\ea
\begin{split}
\phi_{-1}&=-\,\frac{1}{2\sqrt{3}(a+bv+cw)}
(\be_1\overline{y_2y_3\!}\,z_1\!+\!i\ga_1y_1\overline{z_2z_3\!}\,,
\be_2\overline{y_3y_1\!}\,z_2\!+\!i\ga_2y_2\overline{z_3z_1\!}\,,\\
&\qquad\qquad\qquad\qquad\qquad\,\,\,
\be_3\overline{y_1y_2\!}\,z_3\!+\!i\ga_3y_3\overline{z_1z_2\!}\,),
\end{split}
\label{is6eq1}\\
\phi_0&={\ts\frac{1}{\sqrt{3}}}\bigl(y_1z_1,y_2z_2,y_3z_3\bigr),
\label{is6eq2}\\
\begin{split}
\phi_1&={\ts\frac{1}{2\sqrt{3}}}
(\be_1\overline{y_2y_3\!}\,z_1\!-\!i\ga_1y_1\overline{z_2z_3\!}\,,
\be_2\overline{y_3y_1\!}\,z_2\!-\!i\ga_2y_2\overline{z_3z_1\!}\,,\\
&\qquad\quad\,\,\,
\be_3\overline{y_1y_2\!}\,z_3\!-\!i\ga_3y_3\overline{z_1z_2\!}\,).
\end{split}
\label{is6eq3}
\ea
These satisfy
\e
\ms{\phi_{-1}}=(a+bv+cw)^{-1},\quad
\ms{\phi_0}=1 \quad\text{and}\quad
\ms{\phi_1}=a+bv+cw.
\label{is6eq4}
\e

From \eq{is3eq1} and the equation $\ms{\phi_1}=a+bv+cw$ we see that
\begin{equation*}
\phi_2=\frac{\pd\phi_1}{\pd z}
-\frac{\pd}{\pd z}\bigl(\log(a+bv+cw)\bigr)\phi_1.
\end{equation*}
Substituting in for $\phi_1$ from \eq{is6eq3} gives a long and
complicated expression for $\phi_2$. After much calculation using
equations \eq{is5eq2}--\eq{is5eq5}, \eq{is5eq13}--\eq{is5eq14}, 
\eq{is5eq19}--\eq{is5eq21} and other identities satisfied by 
$\be_j,\ga_j$ and $a,b,c$, one can prove that
\e
\phi_2=\xi\phi_{-1}, \quad\text{where}\quad \xi=cC+ibB.
\label{is6eq5}
\e

We can now identify the harmonic sequence of~$\psi$.

\begin{prop} If\/ $bB$ and\/ $cC$ are not both zero then 
$\psi:\R^2\ra\CP^2$ is superconformal, and has harmonic sequence
$(\psi_k)$ given by
\ea
\begin{split}
\psi_{3k-1}(s,t)=
\bigl[&\be_1\overline{y_2y_3\!}\,z_1\!+\!i\ga_1y_1\overline{z_2z_3\!}\,,
\be_2\overline{y_3y_1\!}\,z_2\!+\!i\ga_2y_2\overline{z_3z_1\!}\,,\\
&\be_3\overline{y_1y_2\!}\,z_3\!+\!i\ga_3y_3\overline{z_1z_2\!}\,\bigr],
\end{split}
\label{is6eq6}\\
\psi_{3k}(s,t)=\bigl[&y_1z_1,y_2z_2,y_3z_3\bigr],
\label{is6eq7}\\
\begin{split}
\psi_{3k+1}(s,t)=\bigl[&\be_1\overline{y_2y_3\!}\,z_1\!-\!i\ga_1y_1\overline{z_2z_3\!}\,,
\be_2\overline{y_3y_1\!}\,z_2\!-\!i\ga_2y_2\overline{z_3z_1\!}\,,\\
&\be_3\overline{y_1y_2\!}\,z_3\!-\!i\ga_3y_3\overline{z_1z_2\!}\,\bigr],
\end{split}
\label{is6eq8}
\ea
for all\/ $k\in\Z$. If\/ $bB=cC=0$ then $\psi$ is isotropic, with
finite harmonic sequence $\psi_{-1},\psi_0,\psi_1$ given by
equations \eq{is6eq6}--\eq{is6eq8} with\/~$k=0$.
\label{is6prop1}
\end{prop}

\begin{proof} Since $\phi_2=\xi\phi_{-1}$ where $\xi=cC+ibB$ by
\eq{is6eq5}, if $\xi\ne 0$ then the sequence $(\phi_k)$ exists
for all $k$ and is given by
\begin{equation*}
\phi_{3k-1}=\xi^k\phi_{-1}, \quad
\phi_{3k}=\xi^k\phi_0 \quad\text{and}\quad
\phi_{3k+1}=\xi^k\phi_1.
\end{equation*}
Since $\psi_k=[\phi_k]$, equations \eq{is6eq6}--\eq{is6eq8} follow 
from \eq{is6eq1}--\eq{is6eq3}. Thus $\psi$ is nonisotropic, as
$\psi_k$ exists for all $k$. But any conformal map $\psi:S\ra\CP^2$ is 
isotropic or superconformal from \S\ref{is31}, so $\psi$ is superconformal.

If on the other hand $\xi=0$ then $\phi_2=0$, so $\psi_2$ does not
exist. Thus $\psi$ is isotropic. By \eq{is6eq1}--\eq{is6eq3},
$\psi_{-1},\psi_0$ and $\psi_1$ exist and are given by equations
\eq{is6eq6}--\eq{is6eq8} with $k=0$. But the harmonic sequence
of an isotropic map $\psi:S\ra\CP^m$ has length at most $m+1$,
so this is the whole of the harmonic sequence.
\end{proof}

In the case when $\xi=0$ and $\psi$ is isotropic, $\psi_{-1}$
is holomorphic and $\psi_1$ antiholomorphic. This is not obvious,
but may be proved directly. For instance, when $B=C=0$ we may take the 
$y_j$ and $z_j$ to be {\it real}. Then $\psi$ maps to $\RP^2$ in $\CP^2$, 
and both $\psi_1$ and $\psi_{-1}$ map to the conic $\bigl\{[w_0,w_1,w_2]
\in\CP^2:w_0^2+w_1^2+w_2^2=0\bigr\}$, with~$\psi_{-1}=\overline\psi_1$.

\subsection{Solutions of the Toda lattice and Tzitz\'eica equations} 
\label{is62}

In the rest of the section we assume that $\xi=cC+ibB\ne 0$, so that
$\psi$ is superconformal. Following \S\ref{is32}, we shall construct 
a solution of the Toda lattice equations for $\SU(3)$ out of $\psi$. 
The first thing to do is to find a {\it special}\/ holomorphic
coordinate $z'$ on $\C$, that is, one in which $\xi'=1$ and the $\phi_k'$ 
are periodic with period 3. By \eq{is3eq4}, $z'=z'(z)$ is special if 
\begin{equation*}
\xi'=\bigl({\ts\frac{\pd z'}{\pd z}}\bigr)^{-3}\xi=1.
\end{equation*}
Thus we need $\frac{\pd z'}{\pd z}=\xi^{1/3}$ for some fixed complex 
cube root $\xi^{1/3}$ of $\xi$. So define $z'=\xi^{1/3}(s+it)$. Then
$z'$ is a special holomorphic coordinate on~$\C$.

Working with respect to $z'$ rather than $z$, we get a new sequence
$(\phi_k')$ rather than $(\phi_k)$, with $\phi_k'=-i\xi^{-k/3}\phi_k$.
Thus from \eq{is6eq1}--\eq{is6eq4} we get
\ea
\begin{split}
\phi_{3k\!-\!1}'&=\frac{i\xi^{1/3}}{2\sqrt{3}(a\!+\!bv\!+\!cw)}
\bigl(\be_1\overline{y_2y_3\!}\,z_1\!+\!i\ga_1y_1\overline{z_2z_3\!}\,,
\be_2\overline{y_3y_1\!}\,z_2\!+\!i\ga_2y_2\overline{z_3z_1\!}\,,\\
&\qquad\qquad\qquad\qquad\quad
\be_3\overline{y_1y_2\!}\,z_3\!+\!i\ga_3y_3\overline{z_1z_2\!}\,\bigr),
\end{split}
\label{is6eq9}\\
\phi_{3k}'&=\frac{-i}{\sqrt{3}}\bigl(y_1z_1,y_2z_2,y_3z_3\bigr),
\label{is6eq10}\\
\begin{split}
\phi_{3k\!+\!1}'&=\frac{-i\xi^{-1/3}}{2\sqrt{3}}
\bigl(\be_1\overline{y_2y_3\!}\,z_1\!-\!i\ga_1y_1\overline{z_2z_3\!}\,,
\be_2\overline{y_3y_1\!}\,z_2\!-\!i\ga_2y_2\overline{z_3z_1\!}\,,\\
&\qquad\qquad\quad
\be_3\overline{y_1y_2\!}\,z_3\!-\!i\ga_3y_3\overline{z_1z_2\!}\,\bigr),
\end{split}
\label{is6eq11}\\
\begin{split}
\text{with}&\quad\ms{\phi_{3k-1}'}=\md{\xi}^{2/3}(a+bv+cw)^{-1},\quad 
\ms{\phi_{3k}'}=1 \\
\text{and}&\quad
\ms{\phi_{3k+1}'}=\md{\xi}^{-2/3}(a+bv+cw)
\quad\text{for all $k\in\Z$.}
\end{split}
\label{is6eq12}
\ea

Here we have multiplied by $-i$ because then $\det(\phi_0'\phi_1'
\phi_2')\equiv 1$, as in \eq{is3eq5}. Thus the $\phi_k'$ satisfy 
all the conditions on the $\phi_k$ in \S\ref{is31}--\S\ref{is32}. 
So from \S\ref{is32} if we define $\chi_k=\ms{\phi_k'}$, then the 
$\chi_k$ satisfy the Toda lattice equations for $\SU(3)$ with 
respect to $z'$. Therefore by \eq{is6eq12} we have proved:

\begin{prop} In the situation above, define $\chi_k:\C\ra(0,\iy)$ by
\e
\begin{split}
\chi_{3k-1}&=\md{\xi}^{2/3}(a+bv+cw)^{-1},\quad \chi_{3k}=1
\quad\text{and}\\
\chi_{3k+1}&=\md{\xi}^{-2/3}(a+bv+cw)
\quad\text{for all\/ $k\in\Z$.}
\end{split}
\label{is6eq13}
\e
Then the $\chi_k$ satisfy the Toda lattice equations for $\SU(3)$ 
with respect to $z'=\xi^{1/3}(s+it)$. In terms of\/ $s,t$, this 
means that\/ $\chi_0\chi_1\chi_2\equiv 1$, $\chi_{k+3}=\chi_k$ and
\e
\frac{1}{4\md{\xi}^{2/3}}\left(\frac{\pd^2}{\pd s^2}+
\frac{\pd^2}{\pd t^2}\right)\bigl(\log\chi_k\bigr)=
\chi_{k+1}\chi_k^{-1}-\chi_k\chi_{k-1}^{-1}
\quad\text{for all\/ $k\in\Z$.}
\label{is6eq14}
\e
\label{is6prop2}
\end{prop}

Here \eq{is6eq14} holds because $\frac{\pd^2}{\pd z'\pd\bar z'}=
\frac{1}{4\md{\xi}^{2/3}}\bigl(\frac{\pd^2}{\pd s^2}+
\frac{\pd^2}{\pd t^2}\bigr)$. One can verify \eq{is6eq14} explicitly 
using equations \eq{is5eq13}--\eq{is5eq14}, \eq{is5eq19}--\eq{is5eq21}, 
\eq{is6eq13} and various identities between the $\be_j,\ga_j,B,C,a,b$ 
and $c$. The proposition defines a simple class of doubly-periodic 
solutions $\chi_k$ of the Toda lattice equations for $\SU(3)$. From 
\S\ref{is41} we deduce:

\begin{cor} Define $f:\C\ra(0,\iy)$ by $f=\log(a+bv+cw)-
\frac{2}{3}\log\md{\xi}$. Then $f$ satisfies the Tzitz\'eica 
equation \eq{is4eq2} with respect to~$z'=\xi^{1/3}(s+it)$. 
\label{is6cor}
\end{cor}

Note that the functions $v(s)$, $w(t)$ may be written in terms of 
Jacobi elliptic functions as in \S\ref{is52}, and so the solutions 
in the last two results are entirely explicit. They have a `separated 
variable' form, that is, they are written in terms of single-variable 
functions $v(s)$ and $w(t)$, rather than more general two-variable 
functions $u(s,t)$. The author is not sure whether these solutions 
are already known. 

\subsection{Loops of flat connections and polynomial Killing fields} 
\label{is63}

For the rest of \S\ref{is6} we will work with the special coordinate
$z=\xi^{1/3}(s+it)$, dropping the notation $z'$. From \S\ref{is33}, 
the {\it Toda frame} $F:\R^2\ra\SU(3)$ of $\psi$ is given by 
$F=(f_0f_1f_2)$, where $f_k=\md{\phi_k'}^{-1}\phi_k'$. Using 
equations \eq{is6eq9}--\eq{is6eq12} we may write $F$ down explicitly, 
but we will not do so as the expression is complicated. Then 
$\al=F^{-1}\d F$ is a flat $\SU(3)$ connection matrix on~$\R^2$.

As in \S\ref{is34}, we may extend $\d+\al$ to a loop of flat
$\SU(3)$-connections $\d+\al_\la$ for $\la\in\C$ with $\md{\la}=1$.
We shall write $\al_\la$ out explicitly. Decompose $\al_\la$ as
\e
\al_\la=(\al_1'\la+\al_0')\d z+(\al_{-1}''\la^{-1}+\al_0'')\d\bar z,
\label{is6eq16}
\e
as in \eq{is3eq10}. Then from \eq{is3eq8} and \eq{is6eq13} we find that
\begin{align}
\al_1'&=r^{-1/3}\begin{pmatrix} 0 & 0 & \!\!f^{1/2} \\ 
\!f^{1/2}\!\! & 0 & 0 \\
0 & \!\!rf^{-1}\!\! & 0 \end{pmatrix},\;\> &
\al_0'&=\ha\begin{pmatrix} 0 & 0 & 0 \\ 
0 & \!\!\frac{\pd}{\pd z}(\log f)\!\! & 0 \\
0 & 0 & \!\!-\frac{\pd}{\pd z}(\log f)\! \end{pmatrix},
\label{is6eq17}\\
\al_{\!-1}''&\!=-r^{-1/3}\begin{pmatrix} 0 & \!\!f^{1/2}\!\!\! & 0 \\ 
0 & 0 & \!\!rf^{-1} \\
f^{1/2}\!\!\! & 0 & 0\end{pmatrix},\;\>&
\al_0''&\!=\ha\begin{pmatrix} 0 & 0 & 0 \\ 
0 & \!\!-\frac{\pd}{\pd\bar z}(\log f)\!\! & 0 \\
0 & 0 & \!\!\frac{\pd}{\pd\bar z}(\log f)\! \end{pmatrix},
\label{is6eq18}
\end{align}
where $f=a+bv+cw$ and~$r=\md{\xi}$.

We shall now construct a {\it polynomial Killing field}\/ $\tau$ 
for $\psi$, as in \S\ref{is35}, which is in fact the nontrivial 
polynomial Killing field of lowest degree.

\begin{thm} Write $\xi=r{\rm e}^{i\th}$ for $r>0$ and\/ $\th\in\R$.
Define functions $f,h:\R^2\ra\R$ and\/ $g:\R^2\ra\C$ by
\e
f\!=\!a\!+\!bv\!+\!cw,\;\>
g\!=\!\frac{1}{2f^{1/2}}\Bigl(-b\frac{\d v}{\d s}\!+\!
ic\frac{\d w}{\d t}\Bigr),\;\>
h\!=\!\frac{1}{12f}\Bigl(-b\frac{\d^2 v}{\d s^2}\!
+\!c\frac{\d^2w}{\d t^2}\Bigr),
\label{is6eq19}
\e
and let\/ $\tau=\sum_{n=-2}^2\la^n\tau_n$, where
\begin{align}
\tau_2&=i{\rm e}^{2i\th/3}\begin{pmatrix}
0 & rf^{-1/2} & 0 \\ 0 & 0 & f \\
rf^{-1/2} & 0 & 0 \end{pmatrix},\quad&
\tau_1=i{\rm e}^{i\th/3}\begin{pmatrix}
0 & 0 & g \\ -g & 0 & 0 \\ 0 & 0 & 0 \end{pmatrix}&,
\label{is6eq20}\\
\tau_0&=i\begin{pmatrix}
2h & 0 & 0 \\ 0 & -h & 0 \\ 0 & 0 & -h \end{pmatrix} \quad&
\tau_{-1}=i{\rm e}^{-i\th/3}\begin{pmatrix}
0 & -\bar g & 0 \\ 0 & 0 & 0 \\ \bar g & 0 & 0 \end{pmatrix}&,
\label{is6eq21}\\
&&\mskip -200mu \text{and}\qquad
\tau_{-2}=i{\rm e}^{-2i\th/3}\begin{pmatrix}
0 & 0 & rf^{-1/2} \\ rf^{-1/2} & 0 & 0 \\
0 & f & 0 \end{pmatrix}&.
\label{is6eq22}
\end{align}
Then $\tau$ is a real polynomial Killing field.
\label{is6thm1}
\end{thm}

To prove the theorem one must show that the $\tau_n$ satisfy
\eq{is3eq15} and \eq{is3eq16}. This is a long but straightforward
calculation, using equations \eq{is5eq13}, \eq{is5eq14},
\begin{equation*}
\frac{\pd}{\pd z}=\frac{1}{2r^{1/3}{\rm e}^{i\th/3}}
\Bigl(\frac{\pd}{\pd s}-i\frac{\pd}{\pd t}\Bigr)
\quad\text{and}\quad
\frac{\pd}{\pd\bar z}=\frac{1}{2r^{1/3}{\rm e}^{-i\th/3}}
\Bigl(\frac{\pd}{\pd s}+i\frac{\pd}{\pd t}\Bigr),
\end{equation*}
and identities satisfied by the $\be_j,\ga_j,B,C$ and $\xi$,
and we leave it to the reader.

Both $\al_\la$ and $\tau$ have an extra $\Z_2$-symmetry, which follows
from the fact that $\chi_0\equiv 1$. Define $\ka:\gl(3,\C)\ra\gl(3,\C)$ by
\e
\ka:
\begin{pmatrix} A_{11} & A_{12} & A_{13} \\
A_{21} & A_{22} & A_{23} \\ A_{31} & A_{32} & A_{33} \end{pmatrix}
\mapsto
-\begin{pmatrix} A_{11} & A_{31} & A_{21} \\
A_{13} & A_{33} & A_{23} \\ A_{12} & A_{32} & A_{22} \end{pmatrix}.
\label{is6eq23}
\e
Then $\ka$ is a Lie algebra automorphism, and $\ka^2=1$. It is easy to 
show from \eq{is6eq17}--\eq{is6eq18} and \eq{is6eq20}--\eq{is6eq22} that
\e
\ka(\al_\la)=\al_{-\la} \quad\text{and}\quad 
\ka(\tau(\la))=-\tau(-\la) \quad\text{for all $\la\in\C^*$.}
\label{is6eq24}
\e
The action of $\ka$ on the algebra of polynomial Killing fields
will induce the holomorphic involution $\rho$ on the spectral 
curve discussed in~\S\ref{is42}.

We can now determine the algebra of polynomial Killing fields~$\cal A$.

\begin{thm} In the situation above, the algebra of polynomial Killing 
fields $\cal A$ is generated by $\tau$, $\la^3I$ and $\la^{-3}I$.
\label{is6thm2}
\end{thm}

\begin{proof} Let ${\cal A}'$ be the subalgebra of $\cal A$ generated
by $\tau$, $\la^3I$ and $\la^{-3}I$, and suppose for a contradiction
that ${\cal A}'\ne{\cal A}$. Let $\eta\in{\cal A}\sm{\cal A}'$, and 
take $\eta$ to be real, and of lowest degree $d$. That is, 
$\eta=\sum_{n=-d}^d\la^n\eta_n$ with $\eta_{-n}=-\bar\eta_n^T$
for $n=0,\ldots,d$, and every polynomial Killing field of degree 
less than $d$ lies in~${\cal A}'$.

As $\eta_{d+1}=0$, equations \eq{is3eq15} with $n=d+1$ and
\eq{is3eq16} with $n=d$ show that $\eta_d$ satisfies
\e
[\eta_d,\al_1']=0
\quad\text{and}\quad
\frac{\pd\eta_d}{\pd\bar z}=[\eta_d,\al_0''\,].
\label{is6eq25}
\e
Divide into the three cases (a) $d=3k$, (b) $d=3k+1$, and (c) 
$d=3k+2$ for some $k=0,1,2,\dots$. We will prove a contradiction
in each case in turn.

In case (a), equation \eq{is3eq12} implies that $\eta_d$ is diagonal,
and then as $f$ is nonzero, the first equations of \eq{is6eq17}
and \eq{is6eq25} show that $\eta_d$ is a multiple of the identity.
So write $\eta_d=\ep I$ for some $\ep:\R^2\ra\C$. Taking the trace 
of equations \eq{is3eq15} and \eq{is3eq16} for $n=d$ gives
$\frac{\pd\ep}{\pd z}=\frac{\pd\ep}{\pd\bar z}=0$, as the trace
of any commutator is zero. Thus $\ep$ is constant, and $\eta_d=\ep I$, 
$\eta_{-d}=-\bar\ep I$.

For $k>0$, consider $\eta'=\eta-\ep(\la^3I)^k+\bar\ep(\la^{-3}I)^{-k}$. 
This is a polynomial Killing field of degree less than $d$, as we have
cancelled the terms in $\la^{\pm d}$. Therefore $\eta'\in{\cal A}'$. But 
$\eta=\eta'+\ep(\la^3I)^k-\bar\ep(\la^{-3}I)^{-k}$, so $\eta\in{\cal A}'$,
a contradiction. Also, when $k=0$ we have $\eta=\ep I\in{\cal A}'$. This 
eliminates case~(a).

Similarly, in case (b), equation \eq{is3eq12} and the first equations 
of \eq{is6eq17} and \eq{is6eq25} imply that
\begin{equation*}
\eta_d=\ep r^{-1/3}\begin{pmatrix} 0 & 0 & f^{1/2} \\ f^{1/2} & 0 & 0 \\
0 & rf^{-1} & 0 \end{pmatrix},
\end{equation*}
for some function $\ep:\R^2\ra\C$. The second equation of \eq{is6eq12} is 
equivalent to $\frac{\pd\ep}{\pd\bar z}=0$, so that $\ep$ is holomorphic.
Using the fact that $F_\la\eta F_\la^{-1}$ is independent of $z$ one can 
show that $\ep$ must be constant. This determines $\eta_d$ and~$\eta_{-d}$. 

By \eq{is6eq20}, the leading term of $\tau^2$ is
\begin{equation*}
-\la^4{\rm e}^{4i\th/3}\begin{pmatrix} 0 & 0 & rf^{1/2} \\ 
rf^{1/2} & 0 & 0 \\ 0 & r^2f^{-1} & 0 \end{pmatrix}.
\end{equation*}
Suppose for the moment that $d\ge 7$, so that $k\ge 2$. Consider 
\begin{equation*}
\eta'=\eta+(\la^3I)^{k-1}\xi^{-4/3}\ep\tau^2 
-(\la^{-3}I)^{k-1}\bar\xi^{-4/3}\bar\ep\tau^2.
\end{equation*}
We have cancelled the terms in $\la^{\pm d}$, so $\eta'$ is
a polynomial Killing field of degree less than $d$, and lies 
in ${\cal A}'$. So $\eta$ lies in ${\cal A}'$, a contradiction.

The cases $d=1$ and $d=4$ must be dealt with separately. By
explicit calculation we prove that $\eta$ is a multiple of $I$ 
when $d=1$, and a linear combination of $I,\la^{\pm 3}I$, $\tau$ 
and $\tau^2$ when $d=4$. So $\eta\in{\cal A}'$, finishing case~(b).

In the same way, in case (c) we find that
\begin{equation*}
\eta_d=\ep r^{-2/3}\begin{pmatrix}0 & rf^{-1/2} & 0 \\ 0 & 0 & f \\
rf^{-1/2} & 0 & 0 \end{pmatrix}
\end{equation*}
for some constant $\ep\in\C$. When $d\ge 5$ we define
\begin{equation*}
\eta'=\eta+(\la^3I)^ki\xi^{-2/3}\ep\tau 
+(\la^{-3}I)^ki\bar\xi^{-2/3}\bar\ep\tau,
\end{equation*}
and deduce that $\eta'\in{\cal A}'$, so that $\eta\in{\cal A}'$.
The case $d=2$ we deal with separately, by showing that $\eta$ is 
a linear combination of $\tau$ and $I$, and so lies in ${\cal A}'$.
This completes the proof.
\end{proof}

We can use similar ideas to show that $\psi$ is of {\it finite type},
as in \S\ref{is35}. Define
\e
\eta=(\xi^{-4/3}\la^3-\bar\xi^{-4/3}\la^{-3})\tau^2.
\label{is6eq26}
\e
Then $\eta$ is a real polynomial Killing field of degree 7, and
\eq{is6eq17} and \eq{is6eq20} imply that $\eta_7=\al_1'$ and 
$\eta_6=2\al_0'$. So, by definition, $\psi$ is of finite type.

Furthermore, the proof of the theorem actually implies that
{\it every} polynomial Killing field is of the form
$P_0I+P_1\tau+P_2\tau^2$, where $P_0,P_1,P_2$ are Laurent
polynomials in $\la^{\pm 3}$. Writing $\tau^3$ in this way,
and using the $\Z_2$-symmetry \eq{is6eq24} to eliminate some 
of the terms, we find that $\tau$ must satisfy a cubic equation
\e
\tau^3+D\tau+i(\xi^2\la^6+E+\bar\xi^2\la^{-6})I=0
\label{is6eq27}
\e
for some $D,E\in\R$. Then $\cal A$ is the quotient of the free 
commutative algebra generated by $\la^{\pm 3}I$ and $\tau$ by 
the ideal generated by this equation.

\subsection{The spectral curve}
\label{is64}

Now we can calculate the spectral curve of $\psi$, as in \S\ref{is36}. Define
\begin{equation*}
Y'=\bigl\{(\la,\mu)\in\C^*\t\C:\det\bigl(\mu I-\tau(\la,z)\bigr)=0\bigr\},
\end{equation*}
as in \eq{is3eq19}. Since $\cal A$ is generated by $\tau$ and 
$\la^{\pm 3}I$, this is biholomorphic to the curve  $Y$ of \eq{is3eq17}, 
and so the spectral curve $\ti Y$ as defined by Ferus et al.\ 
\cite[\S 5]{FPPS} is the compactification $\ti Y$ of~$Y'$.

Calculating using \eq{is6eq20}--\eq{is6eq22}, we find that
\ea
\det(\mu I-\tau)&=\mu^3+D\mu+iE+i\xi^2\la^6+
i\bar\xi^2\la^{-6}, \quad\text{where}
\label{is6eq28}\\
D&=f^2+2r^2f^{-1}+2\ms{g}+3h^2\quad\text{and}
\label{is6eq29}\\
E&=-fg^2-f\bar g^2-2f^2h+2r^2f^{-1}h+2\ms{g}h+2h^3.
\label{is6eq30}
\ea
As $Y'$ is independent of $z\in\C$, the functions $D,E$ are constant,
which may be verified directly using \eq{is5eq13}--\eq{is5eq14},
\eq{is6eq19} and identities satisfied by the $\be_j,\ga_j,a,b$,
$c,B,C$ and~$r$.

We can find explicit expressions for these constants by putting $v=w=0$,
which by \eq{is5eq13}--\eq{is5eq14} gives
\begin{equation*}
\bigl({\ts\frac{\d v}{\d s}}\bigr)^2=4(1-B^2),\quad
\bigl({\ts\frac{\d w}{\d t}}\bigr)^2=4(1-C^2)\quad\text{and}\quad
{\ts\frac{\d^2v}{\d s^2}}={\ts\frac{\d^2w}{\d t^2}}=0.
\end{equation*}
Equation \eq{is6eq19} gives values for $f,g$ and $h$, and substituting
these into \eq{is6eq29} and \eq{is6eq30} yields
\e
D=a^2+2a^{-1}(b^2+c^2)\quad\text{and}\quad
E=2\bigl(b^2(1-B^2)-c^2(1-C^2)\bigr).
\label{is6eq31}
\e

This proves that the spectral curve as defined by Ferus et 
al.\ \cite[\S 5]{FPPS} is the compactification $\ti Y$ of 
\e
Y'=\bigl\{(\la,\mu)\in\C^*\t\C:\mu^3+D\mu+iE+i\xi^2\la^6+
i\bar\xi^2\la^{-6}=0\bigr\},
\label{is6eq32}
\e
where $D$ and $E$ are given by \eq{is6eq31}. It can be shown using 
elementary algebraic geometry that $\ti Y$ is nonsingular for
generic $D,E$, with genus 10. Note that the equation satisfied by
$\mu$ in \eq{is6eq32} is the same as that satisfied by $\tau$ 
in~\eq{is6eq27}.

However, McIntosh \cite{McIn1,McIn2,McIn3} uses a different definition 
of the spectral curve. To find it we replace $\la^3$ by $\la$ in 
\eq{is6eq32}, giving 
\e
X'=\bigl\{(\la,\mu)\in\C^*\t\C:\mu^3+D\mu+iE+i\xi^2\la^2+
i\bar\xi^2\la^{-2}=0\bigr\},
\label{is6eq33}
\e
and McIntosh's spectral curve is the compactification $\ti X$ of $X'$.
For generic $D,E$ it is nonsingular with genus 4. The involutions 
$\si:\ti X\ra\ti X$ and $\rho:\ti X\ra\ti X$ discussed in \S\ref{is36} 
and \S\ref{is42} act by
\e
\rho:(\la,\mu)\mapsto(-\la,\mu)\quad\text{and}\quad
\si:(\la,\mu)\mapsto(\bar\la^{-1},-\bar\mu).
\label{is6eq34}
\e

It would be interesting to understand what properties of the 
spectral curve $\ti X$ correspond to the fact that $\psi$ is 
written in terms of single-variable functions $y_k(s)$ and 
$z_k(t)$, rather than more general two-variable functions of 
$(s,t)$. Ian McIntosh has an explanation of this, which may
appear elsewhere.

\subsection{Interpretation using the ideas of \S\ref{is4}}
\label{is65}

Finally we relate the calculations above to the material of 
\S\ref{is4}. From \S\ref{is64} the spectral curve $X$ as defined 
by McIntosh has genus 4. Thus in \S\ref{is43} we have $p=4$ and
$d=2$. The parameter counts there show that the moduli space of 
all finite type genus 4 solutions of the Tzitz\'eica equation,
up to translations in $\R^2$, should have dimension 4. All of
them are expected to be doubly-periodic. For the corresponding
maps $\phi:\R^2\ra{\cal S}^5$ and $\psi:\R^2\ra\CP^2$ to be
doubly-periodic is 4 rationality conditions.

Now the family of genus 4 solutions of the Tzitz\'eica equations
constructed in Corollary \ref{is6cor} depends up to translations
in $\R^2$ on the 3 parameters $\th,B,C$ of \S\ref{is5}. Thus, we 
have not constructed all the genus 4 Tzitz\'eica solutions, but
only a codimension 1 subset of them. This agrees with the analysis 
of \S\ref{is54}, where we were unable to solve the double-periodicity 
conditions in general, because they amounted to 4 rationality 
conditions on 3 variables.

Here are two ways of thinking about why the construction yields only
a codimension 1 subset of the Tzitz\'eica solutions. Firstly, our
solutions have a `separated variable' form, being written in terms
of functions $v(s),w(t)$. It follows that the period vectors of the 
doubly-periodic Tzitz\'eica solutions will point along the $s$ and 
$t$ axes, and so be perpendicular in $\R^2$. However, the general 
genus 4 Tzitz\'eica solution will have period vectors which are not 
orthogonal, and to require them to be orthogonal is a codimension 1 
condition.

Secondly, although the moduli space of quadruples $(\ti X,\rho,\si,\pi)$
with $\ti X$ genus 4 is four-dimensional, the subset which can be defined
by an equation of the form \eq{is6eq33} is only 3-dimensional. In 
\S\ref{is63} we saw that our solutions admit a degree 2 polynomial
Killing field $\tau$, which satisfies a cubic equation over 
$\C[\la^3I,\la^{-3}I]$. It is this cubic equation which gives $X'$
the simple form~\eq{is6eq33}.

So we conclude that although the family of genus 4 Tzitz\'eica 
solutions has dimension 4, only a 3-dimensional subfamily of
these admit a degree 2 polynomial Killing field $\tau$, and it 
is this which is responsible for the special form \eq{is6eq33}
of the spectral curve, and for the other nice behaviour of these
examples. For generic genus 4 Tzitz\'eica solutions the first
non-trivial polynomial Killing field will be of higher degree,
and so the spectral curve will be given by a (singular) equation 
of higher-degree in~$\la^{\pm 2}$.

\section{Extension to three variables}
\label{is7}

Next we generalize Theorem \ref{is5thm1} to a construction of
special Lagrangian 3-folds in $\C^3$ in which all three variables
$r,s,t$ enter in a nontrivial way. The proof is similar to that of 
Theorem \ref{is5thm1}, so we will be brief.

\begin{thm} Let\/ $\al_1,\al_2,\al_3$, $\be_1,\be_2,\be_3$ and\/ 
$\ga_1,\ga_2,\ga_3$ be real numbers with not all\/ $\al_j$, not 
all\/ $\be_j$ and not all\/ $\ga_j$ zero, such that
\e
\begin{aligned}
\al_1\be_1+\al_2\be_2+\al_3\be_3&=0,\quad&
\al_1\ga_1+\al_2\ga_2+\al_3\ga_3&=0,\\
\be_1\ga_1+\be_2\ga_2+\be_3\ga_3&=0
\quad&\text{and}\quad
\al_1\be_1\ga_1+\al_2\be_2\ga_2+\al_3\be_3\ga_3&=0.
\end{aligned}
\label{is7eq1}
\e
Let $I,J,K$ be open intervals in $\R$. Suppose that\/ $x_1,x_2,x_3:I\ra\C$ 
and\/ $u:I\ra\R$ are functions of\/ $r$, that\/ $y_1,y_2,y_3:J\ra\C$ and\/ 
$v:J\ra\R$ are functions of\/ $s$, and\/ $z_1,z_2,z_3:K\ra\C$ and\/ 
$w:K\ra\R$ functions of\/ $t$, satisfying
\begin{align}
\frac{\d x_1}{\d r}&=\al_1\,\overline{x_2x_3\!}\,,&\quad
\frac{\d x_2}{\d r}&=\al_2\,\overline{x_3x_1\!}\,,&\quad
\frac{\d x_3}{\d r}&=\al_3\,\overline{x_1x_2\!}\,,
\label{is7eq2}\\
\frac{\d y_1}{\d s}&=\be_1\,\overline{y_2y_3\!}\,,&\quad
\frac{\d y_2}{\d s}&=\be_2\,\overline{y_3y_1\!}\,,&\quad
\frac{\d y_3}{\d s}&=\be_3\,\overline{y_1y_2\!}\,,
\label{is7eq3}\\
\frac{\d z_1}{\d t}&=\ga_1\,\overline{z_2z_3\!}\,,&\quad
\frac{\d z_2}{\d t}&=\ga_2\,\overline{z_3z_1\!}\,,&\quad
\frac{\d z_3}{\d t}&=\ga_3\,\overline{z_1z_2\!}\,,
\label{is7eq4}\\
\ms{x_1}&=\al_1u+1,&\quad
\ms{x_2}&=\al_2u+1,&\quad
\ms{x_3}&=\al_3u+1,
\label{is7eq5}\\
\ms{y_1}&=\be_1v+1,&\quad
\ms{y_2}&=\be_2v+1,&\quad
\ms{y_3}&=\be_3v+1,
\label{is7eq6}\\
\ms{z_1}&=\ga_1w+1,&\quad
\ms{z_2}&=\ga_2w+1,&\quad
\ms{z_3}&=\ga_3w+1.
\label{is7eq7}
\end{align}
If\/ \eq{is7eq2}--\eq{is7eq4} hold for all\/ $r,s,t$ and\/
\eq{is7eq5}--\eq{is7eq7} hold for some $r,s,t$, then 
\eq{is7eq5}--\eq{is7eq7} hold for all\/ $r,s,t$, for some 
functions $u,v,w$. Define $\Phi:I\t J\t K\ra\C^3$~by
\e
\Phi:(r,s,t)\mapsto\bigl(x_1(r)y_1(s)z_1(t),x_2(r)y_2(s)z_2(t),
x_3(r)y_3(s)z_3(t)\bigr).
\label{is7eq8}
\e
Define a subset\/ $N$ of\/ $\C^3$ by
\e
N=\bigl\{\Phi(r,s,t):\text{$r\in I$, $s\in J$, $t\in K$}\bigr\}.
\label{is7eq9}
\e
Then $N$ is a special Lagrangian $3$-fold in~$\C^3$.
\label{is7thm1}
\end{thm}

\begin{proof} The first part of the theorem, that if 
\eq{is7eq2}--\eq{is7eq4} hold for all $r,s,t$ and \eq{is7eq5}--\eq{is7eq7} 
for some $r,s,t$, then \eq{is7eq5}--\eq{is7eq7} hold for all $r,s,t$, 
follows as in Theorem \ref{is5thm1}. For the second part, we must prove 
that $N$ is special Lagrangian wherever $\Phi$ is an immersion. As in 
Theorem \ref{is5thm1}, this holds if and only if
\ea
\om\Bigl(\frac{\pd\Phi}{\pd r},\frac{\pd\Phi}{\pd s}\Bigr)\equiv
\om\Bigl(\frac{\pd\Phi}{\pd r},\frac{\pd\Phi}{\pd t}\Bigr)\equiv
\om\Bigl(\frac{\pd\Phi}{\pd s},\frac{\pd\Phi}{\pd t}\Bigr)&\equiv 0
\label{is7eq10}\\
\text{and}\qquad \Im\Om\Bigl(\frac{\pd\Phi}{\pd r},
\frac{\pd\Phi}{\pd s},\frac{\pd\Phi}{\pd t}\Bigr)&\equiv 0.
\label{is7eq11}
\ea

Using equations \eq{is7eq2}--\eq{is7eq4} and \eq{is7eq8} we find that
\ea
\frac{\pd\Phi}{\pd r}&=
\bigl(\al_1\,\overline{x_2x_3\!}\,y_1z_1,\al_2\,
\overline{x_3x_1\!}\,y_2z_2,\al_3\,\overline{x_1x_2\!}\,y_3z_3\bigr),
\label{is7eq12}\\
\frac{\pd\Phi}{\pd s}&=
\bigl(\be_1x_1\,\overline{y_2y_3\!}\,z_1,\be_2x_2\,
\overline{y_3y_1\!}\,z_2,\be_3x_3\,\overline{y_1y_2\!}\,z_3\bigr),
\label{is7eq13}\\
\frac{\pd\Phi}{\pd t}&=
\bigl(\ga_1x_1y_1\,\overline{z_2z_3\!}\,,
\ga_2x_2y_2\,\overline{z_3z_1\!}\,,\ga_3x_3y_3\,\overline{z_1z_2\!}\,\bigr).
\label{is7eq14}
\ea
Equations \eq{is7eq12} and \eq{is7eq13} give
\begin{align*}
&\om\bigl({\ts\frac{\pd\Phi}{\pd r},\frac{\pd\Phi}{\pd s}}\bigr)
=\Im(\ov{x_1x_2x_3\!}\,y_1y_2y_3)
\bigl(\al_1\be_1\ms{z_1}+\al_2\be_2\ms{z_2}+\al_3\be_3\ms{z_3}\bigr)\\
&=\Im(\ov{x_1x_2x_3\!}\,y_1y_2y_3)
\bigl(\al_1\be_1(\ga_1w\!+\!1)\!+\!\al_2\be_2(\ga_2w\!+\!1)
\!+\!\al_3\be_3(\ga_3w\!+\!1)\bigr)\\
&=\!\Im(\ov{x_1x_2x_3\!}\,y_1y_2y_3)
\bigl(\al_1\be_1\!+\!\al_2\be_2\!+\!\al_3\be_3\!+\!
w(\al_1\be_1\ga_1\!+\!\al_2\be_2\ga_2\!+\!\al_3\be_3\ga_3)\bigr)\!=\!0,
\end{align*}
using \eq{is7eq7} in the second line and \eq{is7eq1} in the third. This
proves the first equation of \eq{is7eq10}. The second and third follow
in a similar way.

To prove \eq{is7eq11}, observe that
\begin{align*}
\Om\Bigl(\frac{\pd\Phi}{\pd r}&,\frac{\pd\Phi}{\pd s},
\frac{\pd\Phi}{\pd t}\Bigr)=
\Bigl\vert\frac{\pd\Phi}{\pd r}\,\,\frac{\pd\Phi}{\pd s}\,\,
\frac{\pd\Phi}{\pd t}\Bigr\vert=
\left\vert\begin{matrix}
\al_1\,\overline{x_2x_3\!}\,y_1z_1 & \be_1x_1\,\overline{y_2y_3\!}\,z_1 
& \ga_1x_1y_1\,\overline{z_2z_3\!}\, \\
\al_2\,\overline{x_3x_1\!}\,y_2z_2 & \be_2x_2\,\overline{y_3y_1\!}\,z_2 
& \ga_2x_2y_2\,\overline{z_3z_1\!}\, \\
\al_3\,\overline{x_1x_2\!}\,y_3z_3 & \be_3x_3\,\overline{y_1y_2\!}\,z_3 
& \ga_3x_3y_3\,\overline{z_1z_2\!}\,\end{matrix}\right\vert \\
=\bigl(&
\al_1\ms{x_2x_3}\be_2\ms{y_3y_1}\ga_3\ms{z_1z_2}
+\al_2\ms{x_3x_1}\be_3\ms{y_1y_2}\ga_1\ms{z_2z_3}\\
+&\al_3\ms{x_1x_2}\be_1\ms{y_2y_3}\ga_2\ms{z_3z_1}
-\al_1\ms{x_2x_3}\be_3\ms{y_1y_2}\ga_2\ms{z_3z_1}\\
-&\al_2\ms{x_3x_1}\be_1\ms{y_2y_3}\ga_3\ms{z_1z_2}
-\al_3\ms{x_1x_2}\be_2\ms{y_3y_1}\ga_1\ms{z_2z_3}\bigr).
\end{align*}
Thus $\Om\bigl(\frac{\pd\Phi}{\pd r},\frac{\pd\Phi}{\pd s},
\frac{\pd\Phi}{\pd t}\bigr)$ is real, and so~$\Im\Om\bigl(
\frac{\pd\Phi}{\pd r},\frac{\pd\Phi}{\pd s},
\frac{\pd\Phi}{\pd t}\bigr)=0$.
\end{proof}

Here are a few comments on the theorem.
\begin{itemize}
\item[(a)] In Theorem \ref{is5thm1} we took the ranges of $s,t$ 
to be $\R$, but here we take $r,s,t$ in intervals $I,J,K$ in $\R$. 
This is because, by an argument in \cite[Prop.~7.11]{Joyc1}, the 
conditions $\be_1+\be_2+\be_3=0$ and $\ga_1+\ga_2+\ga_3=0$ imply 
that solutions of \eq{is5eq2} and \eq{is5eq3} in some open interval 
extend automatically to all of~$\R$.

However, in Theorem \ref{is7thm1} we do not assume that 
$\al_1+\al_2+\al_3=0$, and so it could happen that $\al_1,\al_2,\al_3$ 
all have the same sign. In this case, solutions $x_j$ to \eq{is7eq2} 
will in general exist in some open interval $I\subset\R$ with 
$\md{x_j}\ra\iy$ at the endpoints of $I$, so that they do not extend 
to $\R$. The same applies to \eq{is7eq3} and~\eq{is7eq4}.

\item[(b)] As in \S\ref{is52} we can write the $x_k,y_k$ and $z_k$ 
entirely explicitly in terms of integrals involving the Jacobi 
elliptic functions. 

\item[(c)] As in Theorem \ref{is5thm2}, in the situation of Theorem 
\ref{is7thm1}, $\frac{\pd\Phi}{\pd r}$, $\frac{\pd\Phi}{\pd s}$ and 
$\frac{\pd\Phi}{\pd t}$ are always complex orthogonal. But in 
general they are not of the same length, so $\Phi$ is not conformal.

\item[(d)] We may recover Theorem \ref{is5thm1} from Theorem 
\ref{is7thm1} as follows. Put $\al_1=\al_2=\al_3=1$, so that 
\eq{is7eq1} becomes equivalent to \eq{is5eq1}. Define 
\begin{equation*}
I=(-\iy,0),\quad x_1(r)=x_2(r)=x_3(r)=-r^{-1} 
\quad\text{and}\quad u(r)=r^{-2}-1,
\end{equation*}
and $J=K=\R$. Then \eq{is7eq2} and \eq{is7eq5} hold, and Theorem 
\ref{is7thm1} becomes equivalent to Theorem \ref{is5thm1}, but 
with a different parametrization for~$r$.
\end{itemize}

\subsection{Description of the family of SL 3-folds}
\label{is71}

We shall now describe the family of SL 3-folds resulting from
Theorem \ref{is7thm1}. We begin by studying the set of solutions 
$\al_j,\be_j,\ga_j$ to \eq{is7eq1}. Define vectors
\begin{align*}
\bs{\al}&=(\al_1,\al_2,\al_3),\quad& 
\bs{\be}&=(\be_1,\be_2,\be_3),\quad&
\bs{\ga}&=(\ga_1,\ga_2,\ga_3),\\
\bs{\al\be}&=(\al_1\be_1,\al_2\be_2,\al_3\be_3),&
\bs{\al\ga}&=(\al_1\ga_1,\al_2\ga_2,\al_3\ga_3),&
\bs{\be\ga}&=(\be_1\ga_1,\be_2\ga_2,\be_3\ga_3)
\end{align*}
in $\R^3$. Rescaling $\bs{\al},\bs{\be}$ and $\bs{\ga}$ has no 
effect on the SL 3-folds constructed in Theorem \ref{is7thm1}, 
so let us assume $\bs{\al},\bs{\be},\bs{\ga}$ are unit vectors. 
We will show that a generic choice of $\bs{\al}$ determines
$\bs{\be},\bs{\ga}$, essentially uniquely.

\begin{prop} Let\/ $\bs{\al}$ be a unit vector in $\R^3$, with\/
$\al_1,\al_2,\al_3$ distinct and nonzero. Then there exist unit
vectors $\bs{\be},\bs{\ga}$ satisfying \eq{is7eq1}, which are
unique up to sign and exchanging~$\bs{\be},\bs{\ga}$.
\label{is7prop}
\end{prop}

\begin{proof} Equation \eq{is7eq1} implies that $\bs{\al},\bs{\be}$ 
and $\bs{\al\be}$ are orthogonal to $\bs{\ga}$. As $\bs{\ga}\ne 0$, 
it follows that $\bs{\al},\bs{\be}$ and $\bs{\al\be}$ are linearly 
dependent. Therefore $\det\bigl(\bs{\al}\,\,\bs{\be}\,\,\bs{\al\be}
\bigr)=0$. This may be rewritten in matrix form as
\e
Q(\bs{\be})\!=\!\ha\!\begin{pmatrix} \be_1 \\ \be_2 \\ \be_3 
\end{pmatrix}^{\!T}\!\!
\begin{pmatrix} 0 & \al_3(\al_2\!-\!\al_1) & \al_2(\al_1\!-\!\al_3) \\
\al_3(\al_2\!-\!\al_1) & 0 & \al_1(\al_3\!-\!\al_2) \\
\al_2(\al_1\!-\!\al_3) & \al_1(\al_3\!-\!\al_2) & 0
\end{pmatrix}
\!\!
\begin{pmatrix} \be_1 \\ \be_2 \\ \be_3 \end{pmatrix}\!\!=\!0.
\label{is7eq15}
\e

Similar equations hold between the $\al_j$ and $\ga_j$, and
between the $\be_j$ and $\ga_j$. Now the $3\t 3$ matrix appearing in 
\eq{is7eq15} has trace zero and determinant $2\al_1\al_2\al_3
(\al_1\!-\!\al_3)(\al_2\!-\!\al_1)(\al_3\!-\!\al_2)$. As by 
assumption $\al_1,\al_2,\al_3$ are distinct and nonzero, this
determinant is nonzero. Hence $Q$ is a trace-free, nondegenerate
quadratic form on~$\R^3$.

Therefore, $\bs{\be}$ must be a unit vector in the intersection of 
the plane $\bs{\al}\cdot\bs{\be}=0$ and the quadric cone \eq{is7eq15} 
in $\R^3$. Let $\bs{\al}^\perp$ be the plane perpendicular to 
$\bs{\al}$, and consider the restriction $Q\vert_{\bs{\al}^\perp}$
of $Q$ to $\bs{\al}^\perp$. As $\bs{\al}$ is a unit vector, we have
\begin{equation*}
0=\Tr(Q)=\Tr\bigl(Q\vert_{\bs{\al}^\perp}\bigr)+Q(\bs{\al}).
\end{equation*}
But $Q(\bs{\al})=0$ by \eq{is7eq15}, so $Q\vert_{\bs{\al}^\perp}$ is 
trace-free.

Thus, by the classification of quadratic forms on $\R^2$, there
exists an orthonormal basis $\bs{\be},\bs{\ga}$ for $\bs{\al}^\perp$ 
such that $Q(x\bs{\be}+y\bs{\ga})=cxy$ for some $c$ and all $x,y$ in
$\R$. If $c=0$ then $Q\vert_{\bs{\al}^\perp}=0$, so $Q$ is degenerate,
a contradiction. So $c\ne 0$, and therefore $\bs{\be},\bs{\ga}$ are
unique up to sign and order, with~$Q(\bs{\be})=Q(\bs{\ga})=0$.

As $\bs{\al},\bs{\be},\bs{\ga}$ are orthonormal they automatically
satisfy the first three equations of \eq{is7eq1}. But by construction
we have arranged that $\bs{\al},\bs{\be}$ and $\bs{\al\be}$ are linearly 
dependent, so $\bs{\al\be}=x\bs{\al}+y\bs{\be}$ for $x,y\in\R$. The
fourth equation of \eq{is7eq1} then follows from the second and third.
\end{proof}

The moral of the proposition is that a generic choice of $\bs{\al}$
determines $\bs{\be}$ and $\bs{\ga}$ up to obvious symmetries.
However, for a nongeneric choice of $\bs{\al}$ there can be more
freedom in $\bs{\be}$ and $\bs{\ga}$. For instance, if we put
$\bs{\al}=3^{-1/2}(1,1,1)$ then $Q\equiv 0$, and $\bs{\be},\bs{\ga}$ 
can be arbitrary orthonormal vectors in~$\bs{\al}^\perp$.

We can now do a parameter count for the family of SL 3-folds coming
from Theorem \ref{is7thm1}. The proposition shows that up to symmetries,
the data $\al_j,\be_j,\ga_j$ has two interesting degrees of freedom.
Also, as in Propositions \ref{is5prop1} and \ref{is5prop2} there exist 
constants $A,B,C\in\R$ such that
\begin{equation*}
\Im(x_1x_2x_3)\equiv A,\quad
\Im(y_1y_2y_3)\equiv B\quad\text{and}\quad
\Im(z_1z_2z_3)\equiv C.
\end{equation*}
Together the $\al_j,\be_j,\ga_j$ and $A,B,C$ determine $N$ up to
automorphisms of $\C^3$. Thus the construction of Theorem \ref{is7thm1}
yields a 5-dimensional family of SL 3-folds, up to automorphisms of~$\C^3$. 

We can also discuss the possible {\it signs} of the $\al_k,\be_k,\ga_k$. 
Suppose for simplicity that $\al_k,\be_k,\ga_k$ are all nonzero. Then 
the four equations of \eq{is7eq1} constrain the signs of $\al_k,\be_k,
\ga_k$, as in each equation the three terms cannot have the same sign,
since their sum is zero. Now permuting $\bs{\al},\bs{\be},\bs{\ga}$, 
and reversing any of their signs, does not change the set of SL 3-folds 
constructed in Theorem~\ref{is7thm1}. 

Considering the constraints on the signs of the $\al_k,\be_k,\ga_k$, 
it is not difficult to show that by permuting and changing signs of 
$\bs{\al},\bs{\be},\bs{\ga}$ we may can arrange that the $\al_k$ are 
all positive, two of the $\be_k$ are positive and one negative, and 
two of the $\ga_k$ positive and one negative. 

With this choice of signs, the argument in \cite[Prop.~7.11]{Joyc1}
shows that solutions $y_k,z_k$ to \eq{is7eq3}--\eq{is7eq4} and
\eq{is7eq6}--\eq{is7eq7} automatically extend to $\R$, so we may take 
$J=K=\R$. However, solutions $x_k$ to \eq{is7eq2} and \eq{is7eq5} 
generally exist only on a proper subinterval $I$ of $\R$. Let us take 
$I$ to be as large as possible.

The discussion of \S\ref{is54} suggests that we should try to arrange
that the $y_k$ are periodic in $s$ and the $z_k$ periodic in $t$. When 
this happens, $\psi$ pushes down to an immersion $I\t T^2\ra\C^3$, whose 
image is a {\it closed}\/ SL 3-fold in $\C^3$. The double-periodicity 
conditions in $s,t$ in this case turn out to be equivalent to those in 
\S\ref{is54}, and there are analogues of parts (a)--(c) of \S\ref{is54} 
in which one can prove they are soluble, which yield countably many 
families of closed, immersed SL 3-folds in $\C^3$ diffeomorphic 
to~$T^2\t\R$.

\subsection{Conclusion: an open problem}
\label{is72}

Theorems \ref{is5thm1} and \ref{is7thm1} are clearly very similar.
But in \S\ref{is6} we saw that the special Lagrangian cones of Theorem
\ref{is5thm1} can be put into a much larger integrable systems
framework. Is there also an `integrable systems' explanation for the 
SL 3-folds of Theorem \ref{is7thm1}? Certainly the solutions of 
Theorem \ref{is7thm1} have many of the hallmarks of integrable systems: 
commuting o.d.e.s, elliptic functions, conserved quantities.

More generally, I suspect that in some sense, SL $m$-folds in $\C^m$ 
for $m\ge 3$ may constitute some kind of higher-dimensional 
integrable system. 

The evidence for this is that there exist many interesting families 
of SL $m$-folds in $\C^m$ which can be written down explicitly, or 
have some other nice properties. For examples, see papers by the author 
\cite{Joyc1,Joyc2,Joyc3,Joyc4,Joyc5}, and others such as Harvey 
and Lawson \cite[III.3]{HaLa}, Haskins \cite{Hask} and Bryant \cite{Brya}.
Also, when the special Lagrangian equations are reduced to an o.d.e., it 
often turns out to be a completely integrable Hamiltonian system, as 
in~\cite[\S 7.6]{Joyc1}.

I have no real idea of how to prove that the special Lagrangian 
equations are integrable, or even of exactly what it would mean for a 
p.d.e.\ to be integrable in more than two dimensions. So I would like 
to bring this question to the attention of the integrable systems 
community, in the hope that someone else may be able to answer it.

\end{document}